\documentclass[ijoc,nonblindrev]{informs3} % current default for manuscript submission

\OneAndAHalfSpacedXII % current default line spacing
\usepackage{natbib}
 \bibpunct[, ]{(}{)}{,}{a}{}{,}%
 \def\bibfont{\small}%
\TheoremsNumberedThrough     % Preferred (Theorem 1, Lemma 1, Theorem 2)
\EquationsNumberedThrough    % Default: (1), (2), ...
\usepackage{booktabs, multirow, array} % for multirow
\usepackage{amsfonts} % for double R
\usepackage{tabularx}
\usepackage{mathtools}
\usepackage{cases}
\usepackage{tikz}
\usetikzlibrary{calc,positioning,shapes.geometric,shapes.symbols,shapes.misc}
\usepackage{mathrsfs}
\usepackage{algorithm}
\usepackage[noend]{algpseudocode}
\usepackage{enumitem}
\usepackage{amssymb}
\usepackage{comment}
\usepackage{csquotes}
\usepackage{indentfirst}
\allowdisplaybreaks[1]

\usepackage{subcaption}
\usepackage{pgfplots}
\usepgfplotslibrary{statistics}
\pgfplotsset{width=10cm,compat=1.9}

\usepackage{lscape}

\begin{document}
%%%%%%%%%%%%%%%%

%\setcounter{page}{1}
%\DOI{10.1287/xxxx.0000.0000}%

\RUNAUTHOR{Jeong et al.}

\RUNTITLE{Adaptive robust electric vehicle routing}

 \TITLE{Adaptive robust electric vehicle routing under energy consumption uncertainty}
%\TITLE{}

\ARTICLEAUTHORS{%
\AUTHOR{Jaehee Jeong}
\AFF{Department of Management Sciences, University of Waterloo, Waterloo, Ontario, Canada, \EMAIL{jaehee.jeong@uwaterloo.ca}}
\AUTHOR{Bissan Ghaddar}
\AFF{Corresponding author. Ivey Business School, Western University, London, Ontario, Canada, \EMAIL{bghaddar@uwaterloo.ca}}
\AUTHOR{Nicolas Zufferey}
\AFF{Geneva School of Economics and Management, University of Geneva, Geneva, Switzerland, \EMAIL{N.Zufferey@unige.ch}}
\AUTHOR{Jatin Nathwani}
\AFF{Department of Management Sciences, University of Waterloo, Waterloo, Ontario, Canada, \EMAIL{nathwani@uwaterloo.ca}}
} % end of the block
\ABSTRACT{%
Electric vehicles (EVs) have been highly favoured as a future transportation mode in the transportation section in recent years. EVs have many advantages compared to traditional transportation, especially the environmental aspect. However, despite many EVs' benefits, operating EVs has limitations in their usage. One of the significant issues is the uncertainty in their driving range. The driving range of EVs is closely related to their energy consumption, which is highly affected by exogenous and endogenous factors. Since those factors are unpredictable, uncertainty in EVs' energy consumption should be considered for efficient operation. This paper proposes an adaptive robust optimization framework for the electric vehicle routing problem. The objective is to minimize the worst-case energy consumption while guaranteeing that services are delivered at the appointed time windows without battery level deficiency. We postulate that EVs can be recharged en route, and the charging amount can be adjusted depending on the circumstance. The proposed problem is formulated as a two-stage adaptive robust problem. A column-and-constraint generation based heuristic algorithm, which is cooperated with variable neighborhood search and alternating direction algorithm, is proposed to solve the proposed model. The computational results show the economic efficiency and robustness of the proposed model, and that there is a tradeoff between the total required energy and the risk of failing to satisfy all customers' demand. 
}%

%\KEYWORDS{deterministic inventory theory; infinite linear programming duality; 
%  existence of optimal policies; semi-Markov decision process; cyclic schedule}


\maketitle
%%%%%%%%%%%%%%%%%%%%%%%%%%%%%%%%%%%%%%%%%%%%%%%%%%%%%%%%%%%%%%%%%%%%%%

% Samples of sectioning (and labeling) in IJOC
% NOTE: (1) \section and \subsection do NOT end with a period
%       (2) \subsubsection and lower need end punctuation
%       (3) capitalization is as shown (title style).
%
%\section{Introduction.}\label{intro} %%1.
%\subsection{Duality and the Classical EOQ Problem.}\label{class-EOQ} %% 1.1.
%\subsection{Outline.}\label{outline1} %% 1.2.
%\subsubsection{Cyclic Schedules for the General Deterministic SMDP.}
%  \label{cyclic-schedules} %% 1.2.1
%\section{Problem Description.}\label{problemdescription} %% 2.

% Text of your paper here

\section{Introduction}
Over the last two decades, the transportation sector has rapidly grown. For instance, the number of the registered vehicles in US has increased by over 43\% in the last 20 years (from 1990 to 2019) (\cite{statista}). However, this rapid growth also brings some adverse effects, especially environmental problems. United States Environmental Protection Agency (US EPA) reported that the transportation sector accounts for the most significant portion (29\%) of the total greenhouse gas (GHG) emissions of US in 2019 (\cite{epa2019}). Also, \cite{wri2021} reported that the transportation sector is the one of the fastest-growing sources of GHG emissions, which has grown by 79\% from 1990, and it stands out 16.2\% of total global GHG emissions in 2019. Consequently, the interest in green vehicle routing has increased to consider environmental challenges facing the transportation sector.

In the recent decade, electric vehicles (EVs) have become the abiding interest of the transportation sector as a future mode of transportation because of their various benefits: reducing oil dependency, carbon dioxide (CO$_2$) emission reduction, less noise, and utilization of renewable energy (\cite{bradley2009design, sioshansi2009emissions, sundstrom2010optimization, yi2017effects}). Indeed, the number of EVs on the world's roads increased from 17,000 in 2010 to more than 10 million in 2020, and its number is rapidly growing worldwide (\cite{iea2021techreport}).  Additionally, many leading companies, such as UPS, FedEx, and Walmart, have deployed EV fleets in their operations (\cite{winstoninside}). However, despite several benefits of EVs and their growing penetration, EVs still have critical issues limiting their usage. One issue is their limited driving range compared to conventional vehicles fueled by gas or diesel. %EVs are powered by a battery, which is typically big and heavy and because of the limited space of the cabin, only a limited number of batteries can be loaded, which leads to the limited driving range. 
Other problems include the scarcity of the charging infrastructure and the long recharging time. Therefore, a more careful routing plan is needed to overcome EV's limited driving range, especially for commercial EVs.

The driving range of EVs depends on the vehicle's battery capacity and its efficiency. Many factors relate to the vehicle's energy consumption efficiency, such as vehicle-, environmental- and driver-related factors, and several prediction models were proposed to predict EV's energy consumption depending on those factors (\cite{zhang2020energy, basso2021electric}). The developed forecasting models could provide good predictions to the EV's energy consumption, but an exact forecast of the EV's energy consumption is not possible since environmental factors have a high degree of uncertainty. One of the major factors that affect EV's energy consumption efficiency is temperature. The temperature directly impacts the battery performance, but it also relates to auxiliary heating and cooling. Cold or hot days make people use the air conditioning to heat or cool the vehicle. The net effects of temperature and the resulted auxiliary power demand can decrease EV's driving range up to 40\% compared to its maximum potential range (\cite{yuksel2015effects}). %Especially, the impacts of auxiliary energy demand have a higher impact during low-speed driving ($\le$ 80km) and long-time travel (\cite{asamer2016sensitivity}). 
%Weather conditions are affecting the EV range, especially the one from the on-board climate control systems depending on environmental conditions like high or low outside temperatures (AC cooling or heating), precipitation (window wipers), visibility (car lights), snowing and freezing conditions (windshield or park heaters). It is important is to have a frequently updated weather forecast model (e.g., GFS) with very high spatial and temporal resolution to access forecast data and to directly feed these, together with local historical data, into predictive analytics, with the aim of extracting the weather effect on the EV energy consumption.
Furthermore, even if the temperature can be precisely predicted, the EV's performance could be 32\% higher or less than the expected driving range of that temperature depending on other factors, such as wind speed, rolling resistance, terrain, driver habits, trip length and start conditions (\cite{yi2017effects, geotab}). Therefore, to operate EVs efficiently, the uncertainty in EVs' energy consumption should be considered.

The purpose of this paper is to formulate and solve a model for a practical delivery service system in which EVs are used for commercial vehicles delivering products to customers, and all routes are guaranteed to withstand any realization of energy consumption uncertainties. To deal with uncertainty, we adopt an adaptive robust optimization framework (\cite{ben2004adjustable}). Thus, the proposed model is an adaptive robust electric vehicle routing problem with a time window, partial recharge, and energy consumption rate uncertainty (AR-EVRPTWPR). It extends the electric vehicle routing problem with a time window and partial recharge (EVRPTWPR) proposed by \cite{keskin2016partial} by considering the energy consumption rate uncertainty. The AR-EVRPTWPR is formulated as a two-stage adaptive robust problem and solved by means of a column-and-constraint generation framework (CCG) {(\cite{zeng2013solving, Yu2021})}. Since the CCG decomposes the model into a master problem and subproblem, we propose a hybrid Variable Neighborhood Search--Tabu Search (VNS-TS) metaheuristic (\cite{schneider2014electric}) and Alternating Direction (AD) algorithm (\cite{konno1976cutting}) to solve the master problem and subproblem, respectively. The proposed method is applied to test instances of \cite{schneider2014electric} to show the performance of the approach. 

%This paper is the first to apply adaptive robust optimization to EVRPTWPR with energy consumption uncertainties. The energy consumption uncertainty has a significant impact on the operation that uses EVs as the primary vehicle. The proposed model provides

The remainder of this paper is organized as follows. Section \ref{literature_review} presents a literature review of related work. Section \ref{Formulations} illustrates an example of the limitation of robust optimization for EVRPTWPR and provides the mathematical formulation of AR-EVRPTWPR. The proposed solution methods are presented in Section \ref{solutionmethod}. In Section \ref{results}, computational experiments are conducted for the proposed AR-EVRPTWPR. Finally, conclusions and future research directions are drawn in Section \ref{conclusion}.

\section{Literature Review} \label{literature_review}

Since 1959, when \cite{dantzig1959truck} first proposed the vehicle routing problem (VRP), VRP has been extensively studied for its variants and solution methods (\cite{eksioglu2009vehicle, toth2014vehicle, braekers2016vehicle}). The objective of VRP is to find a set of routes for a fleet of vehicles that minimizes the total travel distances while it starts from the depot, visits a given set of customers, and returns to the depot. Recently, the green vehicle routing problem, a variation of VRP, has been widely studied due to the negative effects of the transportation on the environment. A comprehensive review of green vehicle routing problems can be found in \cite{demir2014review, lin2014survey, bektacs2019role}.

The electric vehicle routing problem (EVRP) is a variation of the VRP where EVs are integrated into the distribution operations (\cite{erdougan2012green}). Many variants of EVRP have been studied and a recent solver was developed by \cite{Kullman2021}. \cite{conrad2011recharging} first proposed a recharging VRP allowing vehicle recharging en tour. \cite{schneider2014electric} considered time windows in EVRP and where the EVs' batteries are only allowed to fully recharge at charging stations. This EVRP with time windows and recharge stations was extended to consider partial recharge by \cite{felipe2014heuristic}, \cite{keskin2016partial}, and \cite{schiffer2017electric}. \cite{hiermann2019routing} considered a fleet mix with several types of vehicles. \cite{cortes2019electric} regarded the recharging time as an idle time and allowed for servicing the customers by walking from the charging station while the EV is recharging. 

The above studies consider deterministic parameters, i.e., the expected or average value of energy consumption, travel time, demand, and recharging time. Thus, a solution obtained from a deterministic EVRP would be very vulnerable to variations in the parameters. Some stochastic variants of EVRP have been proposed to deal with the uncertainty in EVRP, which minimized the expected objective under demand, travel time or service time uncertainty (\cite{gendreau201650th, keskin2021simulation}). Also, a robust optimization model is proposed by \cite{pelletier2019electric} to take into account energy consumption uncertainty in EVRP with time windows. \cite{basso2021electric} applied a probabilistic Bayesian machine learning approach to predict energy consumption and modelled the problem as an EVRP with chance-constrained partial recharging problem.

Adaptive (or adjustable) robust optimization (ARO) is a branch of robust optimization (RO), which recently has an increased range of applications (\cite{ben2004adjustable}). ARO allows the subset of decision variables to adapt itself against uncertainty scenarios in the uncertainty set. The comprehensive studies of ARO approach can be found in application of the power system operation (\cite{bertsimas2012adaptive, wang2013two, lorca2014adaptive}). Also, there exist ARO approaches for conventional vehicle routing problems (\cite{lee2012robust, agra2013robust}) and inventory routing problems (\cite{agra2018robust}). To the best of our knowledge, there has not been any research considering energy consumption uncertainty in EVRP using an ARO approach. In the context of EVRP, energy consumption uncertainty was considered in \cite{pelletier2019electric} and \cite{basso2021electric}, but the former adopts the RO approach that only considers the worst energy consumption and the latter accounts for chance constraints with energy consumption prediction through machine learning techniques. The RO approach usually provides a conservative solution, and the solution quality of the chance constraints approach highly depends on which distribution function is used for. However, ARO generally provides a less conservative solution than the classical RO approach because ARO separates decision variables into two groups: here-and-now and wait-and-see decisions. Here-and-now decisions should be determined before the uncertainty realization, but wait-and-see decisions can adapt to uncertainty realization. This ARO's characteristic is similar to the nature of the EV routing problem under consideration. Each EV's route should be decided at least a day ahead, but the refuelling amount can be decided based on the operating situation. Therefore, ARO approach is more suitable to consider uncertainties in routing problems. Also, the solution quality of the ARO approach does not depend on the distribution function like the chance-constrained approach. Instead, the ARO approach requires more simple statistics, such as mean and variances. Thus, the ARO approach does not depend on estimating the precise distribution which might not be possible in some applications.

The main contributions of this paper are to introduce, model, and solve an adaptive robust optimization model of a delivery routing problem where a fleet of EVs need to visit a given set of customers regardless of any realization of energy consumption scenario within the uncertainty set. In the next section, we illustrate a simple example of the limitation of the classical robust optimization approach to tackle this problem. Our model extends the work of CCG to apply the adaptive robust optimization approach. The proposed model is formulated as a two-stage adaptive robust model, and a decomposition solution framework is proposed to solve our model.

\section{EVRPTWPR Formulations} \label{Formulations}
In this section, we describe the motivation for developing an adaptive robust EVRPTWPR with energy consumption rate uncertainty and then present its mathematical formulation. In the following sections, three models related to the EVRPTWPR model are presented: the deterministic EVRPTWPR model, the robust EVRPTWPR with energy consumption rate uncertainty (RO-EVRPTWPR) model, and the proposed adaptive robust EVRPTWPR with energy consumption rate uncertainty (AR-EVRPTWPR) model. 

\subsection{Deterministic EVRPTWPR Formulation}
In this section, we first start by presenting the deterministic EVRPTWPR model. We use the same deterministic EVRPTWPR model as the one presented in \cite{keskin2016partial} and we follow a similar notation (as shown in Table \ref{notation}). The deterministic model with time window constraints and partial recharging is formulated as follows:
\begin{subequations}
\begin{align}
    \min ~~~& \sum_{i \in V_0^{'}, j \in V_{N+1}^{'}, i \neq j} d_{ij} x_{ij}
    \label{keskinobj}\\
    \textrm{s.t.} ~~~& \sum_{j \in V_{N+1}^{'}, i \neq j} x_{ij} = 1 && \forall i \in V \label{keskinconb}\\
    & \sum_{j \in V_{N+1}^{'}, i \neq j} x_{ij} \le 1 && \forall i \in F^{'} \label{keskinconc}\\
    & \sum_{i \in V_0^{'}, i \neq j} x_{ij} - \sum_{i \in V_{N+1}^{'}, i \neq j} x_{ji} = 0 && \forall j \in V^{'} \label{keskincond}\\
    & \tau_i + (t_{ij} + s_i)x_{ij} - l_0(1- x_{ij}) \le \tau_j && \forall i \in V_0, j \in V_{N+1}^{'}, i \neq j \label{keskincone}\\
    & \tau_i + t_{ij}x_{ij} + g(Y_i - y_i) - (l_0 + gQ)(1 - x_{ij}) \le \tau_j && \forall i \in F^{'}, j \in V_{N+1}^{'}, i \neq j \label{keskinconf}\\
    & e_j \le \tau_j \le l_j && \forall j \in V_{0, N+1}^{'} \label{keskincong}\\
    & 0 \le u_j \le u_i - q_i x_{ij} + C(1 - x_{ij}) && \forall i \in V_{0}^{'}, j \in V_{N+1}^{'}, i \neq j \label{keskinconh}\\
    & 0 \le u_0 \le C &&\label{keskinconi}\\
    & 0 \le y_j \le y_i - (h \cdot d_{ij})x_{ij} + Q(1 - x_{ij}) && \forall i \in V, j \in V_{N+1}^{'}, i \neq j \label{keskinconj}\\
    & 0 \le y_j \le Y_i - (h \cdot d_{ij})x_{ij} + Q(1 - x_{ij}) && \forall i \in F_0^{'}, j \in V_{N+1}^{'}, i \neq j \label{keskinconk}\\
    & y_i \le Y_i \le Q && \forall i \in F_0^{'} \label{keskinconl}\\
    & x_{ij} \in \{0,1\} && \forall i \in V_0^{'}, j \in V_{N+1}^{'}, i \neq j, \label{keskinconm}
\end{align}
\label{keskinmodel}
\end{subequations}
where $V, F$, and $F^{'}$ denote the sets of customers, recharge stations, and dummy recharge stations, respectively, and $V^{'} := V \cup F^{'}, V_0^{'} := V^{'} \cup \{0\}, V_{N+1}^{'} := V^{'} \cup \{N+1\}$, and $F_0^{'} := F^{'} \cup \{0\}$ where $0$ and $N+1$ denote the initial depot and the end depot, respectively. Decision variable $x_{ij}$ is the binary decision of travel of arc $(i,j)$, $u_j$ is the cargo level decision at vertex $j$, $\tau_j$ and $y_j$ specify the service start time and battery level when the vehicle arrives at vertex $j$, respectively, and $Y_j$ is the battery level when the vehicle departs from vertex $j$. Parameters $C, Q, g$, and $h$ are the cargo capacity, battery capacity, recharging rate, and energy consumption rate, respectively. Also, $e_j$ and $l_j$ are the earliest and latest service start time at vertex $j$. The details of sets, variables, and parameters used in this paper are summarized in Table \ref{notation}. The objective (\ref{keskinobj}) minimizes the travel distance. Constraints (\ref{keskinconb}) correspond to the connectivity of customer visits. Constraints (\ref{keskinconc}) enforce the connectivity of visits to recharging station. Constraints (\ref{keskincond}) correspond to the flow conservation for each vertex, where the number of incoming and outgoing arcs should be the same. Constraints (\ref{keskincone}) and (\ref{keskinconf}) enforce arcs' time feasibility leaving from vertex in $V_0$ and $F^{'}$, respectively. Constraint (\ref{keskincong}) ensure the time windows of each vertex. Constraints (\ref{keskinconh}) and (\ref{keskinconi}) enforce the demand fulfillment of all customers, and cargo level should be nonnegative and cannot exceed the vehicle's capacity. Constraints (\ref{keskinconj}) and (\ref{keskinconk}) enforce that the vehicle's battery charge level is always nonnegative. Constraints (\ref{keskinconl}) state that the vehicle's battery charge level departing vertex $j$ should be larger than its battery charge level when it arrives and cannot exceed its battery capacity.

\begin{table}
\caption{Summary of notation}
    \centering
    \small
    \begin{tabularx}{\textwidth}{lX}
    \toprule
    \textbf{Sets} \\
    $0, N+1$ & Initial depot and end depot \\
    $V$     & Set of customer vertices $V := \{1,2,...,N\}$ \\
    $F$     & Set of recharge station vertices \\
    $F^{'}$     & Set of recharge station vertices and dummy vertices of the set of recharging station $F$\\
    $V^{'}$     & Set of customer, recharge station, and dummy recharge station vertices $V^{'} := V \cup F^{'}$ \\
    $V_0^{'}$ & $V_0^{'} := V^{'} \cup \{0\}$ \\
    $V_{N+1}^{'}$ & $V_{N+1}^{'} := V^{'} \cup \{N+1\}$ \\
    $V_{0, N+1}^{'}$ & $V_{0, N+1}^{'} := V^{'} \cup \{0, N+1\}$ \\
    $F_0^{'}$ & $F_0^{'} := F^{'} \cup \{0\}$ \\
    $A$ &   Set of arcs $A := \{(i,j) ~ | ~ i \in V_0^{'}, j \in V_{N+1}^{'}, i \neq j\}$ \\
    $A(\rho)$ & Set of arcs in the route $\rho$ \\
    $\mathcal{U}$ & Set of uncertain energy consumption rate \\
    \\
    \textbf{Variables} \\
    $x_{ij}$    &   Binary decision variable for the arc $(i,j)$ travelling.\\
    $u_j$    &   Nonnegative continuous decision variable for remaining cargo level when the vehicle arrives the vertex $j$\\
    $\tau_j$    &   Continuous decision variable for the service start time of the vertex $j$\\
    $y_j$    &   Nonnegative continuous decision variable for battery charge level when the vehicle arrives the vertex $j$\\
    $Y_j$    &   Nonnegative continuous decision variable for battery charge level when the vehicle departs the vertex $j$\\
    $p_j$    &   Nonnegative continuous decision variable for batter recharging amount at the vertex $j$\\
    \\
    \textbf{Parameters} \\
    $d_{ij}$    &   Distance of arc $(i,j)$\\
    $t_{ij}$    &   Travel time of arc $(i,j)$\\
    $q_j$   &   Demand of vertex $j$\\
    $e_j, l_j$   &   Earliest and latest service start time of vertex $j$\\
    $C$   &   Vehicle capacity\\
    $Q$   &   Vehicle battery capacity\\
    $g$   &   Recharging rate\\
    $h$   &   Energy consumption rate\\
    \bottomrule
    \end{tabularx}
    \label{notation}
\end{table}

\subsection{RO-EVRPTWPR Formulation and Limitations} \label{ROEVRPTWPR}
We note that up to our knowledge, there is no research that presents the RO-EVRPTWPR formulation considering energy consumption rate uncertainty and thus we provide the formulation in this section. Furthermore, we show a crucial limitation of the classical robust optimization approach which motivates the development of an adaptive robust optimization approach for the EVRPTWPR. The classical robust optimization approach deals with the uncertainty of data by finding a \textit{robust solution} that is immunized against all possible variations in data (\cite{lee2012robust}). Here, the \enquote{\textit{immunized}} denotes that a solution should be feasible to all possible scenarios, and decision variables are not allowed to adjust for the realized uncertainty. Unfortunately, this nature of classical robust optimization easily causes the solution to be conservative, or it sometimes fails to find a robust feasible solution. The details and extensive literature review of the robust optimization can be found in for example in \cite{ben2004adjustable} and \cite{bertsimas2011theory}.

To consider energy consumption rate uncertainty in EVRPTWPR, every constraint related to the energy consumption rate becomes a robust constraint. Thus, constraints (\ref{keskinconj}) and (\ref{keskinconk}) should be robust against energy consumption rate uncertainty. Robust constraints can be rewritten as follows:
\begin{subequations}
\begin{align}
    & 0 \le y_j \le y_i - (h_{ij} \cdot d_{ij})x_{ij} + Q(1 - x_{ij}) && \forall i \in V, j \in V_{N+1}^{'}, i \neq j, \forall \boldsymbol{h} \in \mathcal{U} \label{robustcona}\\
    & 0 \le y_j \le Y_i - (h_{ij} \cdot d_{ij})x_{ij} + Q(1 - x_{ij}) && \forall i \in F_0^{'}, j \in V_{N+1}^{'}, i \neq j, \forall \boldsymbol{h} \in \mathcal{U}, \label{robustconb}
\end{align}
\label{robustcon}
\end{subequations}
where $\boldsymbol{h}$ and $\mathcal{U}$ are the vector of uncertain energy consumption rate and the set of uncertain parameters. Since we consider energy consumption uncertainty, each arc has its own energy consumption rate $h_{ij}$. Constraints (\ref{robustcona}) and (\ref{robustconb}) enforce that the vehicle's battery charge level is non-negative for all possible realizations of $\boldsymbol{h}$. We note that variation in energy consumption directly affects to vehicle's battery charge level and the affected battery charge level in turn affects the recharging time at the charging stations. In addition, the change of recharging time may impact the service start time. Therefore, the above robust constraints implicitly affect the other constraints, such as constraints (\ref{keskincone}) -- (\ref{keskincong}), and (\ref{keskinconl}). To illustrate that the classical robust optimization approach is not suitable to the EVRPTWPR, we next describe the following simple example.

\begin{figure}[!ht]
    \centering
    \begin{tikzpicture}[roundnode/.style={circle, draw=green!60, fill=green!5, thick, minimum size=3mm}, squarednode/.style={rectangle, draw=red!60, fill=red!5, thick, minimum size=3mm}, trianglenode/.style={isosceles triangle,isosceles triangle apex angle=60, shape border rotate=90, draw=blue!60, fill=blue!5, thick, minimum size=1mm}]
        %Nodes
        \node[squarednode, scale = 1]      (depot)                              {\footnotesize D};
        \node[roundnode, scale = 0.85]        (c1)       [right =5cm of depot] {\footnotesize C1};
        \node[roundnode, scale = 0.85]      (c2)       [below=3cm of depot] {\footnotesize C2};
        \node[trianglenode, scale = 0.5]        (s1)       [below right=1.694cm of depot] {\Large S1};

        %Lines
        \draw[-] (c1.west) -- (depot.east) node[above, midway] {\footnotesize 40};
        \draw[-] (depot.south) -- (c2.north) node[left, midway] {\footnotesize 30};
        \draw[-] (c1.south west) -- (c2.north east) node[below=1mm, midway] {\footnotesize 50};
        \draw[-] (depot.south east) .. controls (1.0,-0.8).. (s1.north west) node[right, midway] {\footnotesize 19};
        \draw[-] (c1.200) .. controls (3,-1.4).. (s1.35) node[above left=-1mm, midway] {\footnotesize 35};
        \draw[-] (c2.65) .. controls (0.4,-2.5).. (s1.215) node[above=1mm, midway] {\footnotesize 20};
        % \draw[->] (rightsquare.south) .. controls +(down:7mm) and +(right:7mm) .. (lowercircle.east);
        
        \node[above=0mm of depot] (time) {\scriptsize$[0,140]$};
        \node[above=-0.5mm of c1] (time) {\scriptsize$[40,50]$};
        \node[below=0mm of c2] (time) {\scriptsize$[90,100]$};
        \node[below=0mm of s1] (time) {\scriptsize$[0,140]$};
\end{tikzpicture}
\caption{An example of two customers and one charging station}
    \label{figure_example_1}
\end{figure}
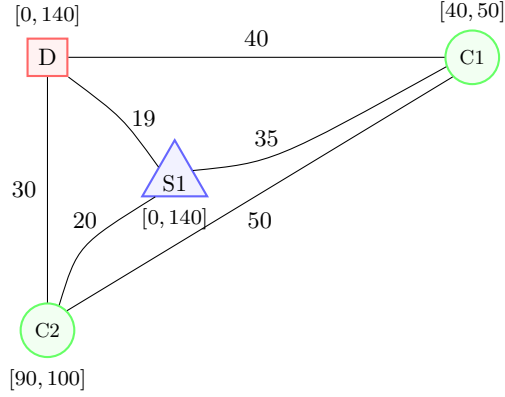

\begin{example}
    Consider a small network that consists of two customer vertices and one charging station. Initial and end depot are the same vertex, and two customers need to be served. The network is depicted in Figure \ref{figure_example_1}. The travel times are indicated on each arc as shown in Figure \ref{figure_example_1}, and the time windows are also denoted on each vertex. For ease of explanation, we assume a service time of zero for all vertices. Let the battery capacity $Q = 120$, and the recharging rate be $1.0$. The expected energy consumption rate is assumed to be 1.0 for all arcs, and we consider six energy consumption rate scenarios. In each scenario, a single arc has a 10\% more energy consumption rate than the average value. If any of $(D, C1), (C1, C2)$, or $(C2, D)$ arc has an energy consumption rate of 1.1, the vehicle should visit the recharge station $S1$. The shortest route visiting $S1$ is $R^{0} = D \rightarrow C1 \rightarrow S1 \rightarrow C2 \rightarrow D$, but $R^0$ is infeasible when $h_{D, C1} = 1.1$. For instance, $y_{C1}$ and $y_{S1}$ are $76$ and $41$, respectively, and the recharging time is at least $9$ ($50$ needs to travel arcs $(S1, C2)$ and $(C2, D)$ and $50 - 41 = 9$). The earliest arrival time to $C2$ $\tau_{C2}$ is $104$, which violates the $C2$'s time window constraint. Thus, to be feasible for both scenarios $h_{D, C1} = 1.1$ and $h_{C1, C2} = 1.1$, the vehicle should travel through route $R^1 = D \rightarrow C1 \rightarrow C2 \rightarrow S1 \rightarrow D$, and $\boldsymbol{y}^{1} = (y_{C1}, y_{C2}, y_{S1}, Y_D, Y_{S1})$ and $\boldsymbol{\tau}^{1} = (\tau_{C1}, \tau_{C2}, \tau_{S1})$ have a value of $(76, 21, 1, 120, 19)$ and $(40, 90, 109)$, respectively. However, $y^1$ becomes infeasible under a scenario of $h_{C2, S1} = 1.1$. For instance, $y_{S1} \le y_{C2} - 22 = 21 - 22$, so $y_{S1} \le -1$ means that the vehicle arrives at $S1$ with a negative battery level. There is no more feasible route. This result shows that RO-EVRPTWPR has no feasible solution for Figure \ref{figure_example_1}.
    
    %Old one 
    %Then, we can easily infer that the worst case scenario is a scenario with $h_{C1, C2} = 1.1$. The optimal route $R^{1*}$ will be $D \rightarrow C1 \rightarrow S1 \rightarrow C2 \rightarrow D$. and the optimal values of $\boldsymbol{y}^{1*} = (y_{C1}, y_{C2}, y_{S1}, Y_D, Y_{S1})$ and $\boldsymbol{\tau}^{1*} = (\tau_{C1}, \tau_{C2}, \tau_{S1})$ are $(80, 30, 45, 120, 49)$ and $(40, 94, 75)$, respectively. However, $\boldsymbol{y}^{1*}$ violates constraint (\ref{robustconb}) when $h_{D, C1} = 1.1$, for instance, $80 \le 120 - 44$ does not hold. To be feasible for both scenarios of $h_{C1, C2} = 1.1$ and $h_{D, C1}=1.1$, $y_{C1}$ and $y_{S1}$ should decrease to 76 and 41, respectively, and then recharging time at $S1$ increases to 8 (i.e., $Y_{S1} - y_{S1} = 49 - 41$). As a result, $\tau_{C1}$ becomes 102, and this leads to the infeasibility.
\end{example}

\noindent From a mathematical model point of view, we can claim the following result:
\begin{claim}
    The classical robust optimization model for the EVRPTWPR with energy consumption rate uncertainty is the same as that every arc has its worst energy consumption realization.
    \proof{Proof}
        For arbitrary arc $(i,j) \in A$, let consider the constraint (\ref{robustcona}) when the uncertainty set $\mathcal{U} = \{\boldsymbol{h}^1,...,\boldsymbol{h}^n\}$. Then, we can explicitly write the constraint (\ref{robustcona}) as follows:
        \begin{align}
        \begin{split}
            & 0 \le y_j \le y_i - (h^1_{ij} \cdot d_{ij})x_{ij} + Q(1 - x_{ij}) \\
            & \hspace{3cm} \vdots \\
            & 0 \le y_j \le y_i - (h^n_{ij} \cdot d_{ij})x_{ij} + Q(1 - x_{ij}) 
        \end{split}
        \label{claim_1_pf_1}
        \end{align}
        Since every constraint is redundant when $x_{ij} = 0$, we only need to consider a case of $x_{ij} = 1$. Suppose $x_{ij} = 1$, a constraint with $h^{max}_{ij} = \max_{k=1,...,n} h^k_{ij}$ dominates the other constraints. It also holds for the constraint (\ref{robustconb}). Thus, constraints (\ref{robustcon}) can be rewritten as follow:
        \begin{subequations}
        \begin{align}
            & 0 \le y_j \le y_i - \bigg( \max_{k=1,...,n} h^k_{ij} \cdot d_{ij} \bigg)x_{ij} + Q(1 - x_{ij}) && \forall i \in V, j \in V_{N+1}^{'}, i \neq j \label{robustcona_alt}\\
            & 0 \le y_j \le Y_i - \bigg(\max_{k=1,...,n} h^k_{ij} \cdot d_{ij}\bigg)x_{ij} + Q(1 - x_{ij}) && \forall i \in F_0^{'}, j \in V_{N+1}^{'}, i \neq j \label{robustconb_alt}
        \end{align}
        \label{robustcon_alt}
        \end{subequations}
        \hfill\Halmos
    \endproof
\end{claim}

This result shows that the classical robust optimization approach is not suitable for the EVRPTWPR with energy consumption uncertainty. Adopting the classical robust optimization approach can lead to a conservative optimal solution or fail to find a feasible solution. To properly handle the uncertainty in EVRPTWPR, some decision variables should be allowed to vary to deal with the different uncertain realizations. It motivates us to consider the AR-EVRPTWPR as described in the next section.

\subsection{AR-EVRPTWPR Formulation}

In this section, we propose an Adaptive Robust Optimization (ARO) model for EVRPTWPR with energy consumption rate uncertainty. Instead of using a two-index formulation, we incorporate the index of vehicle $k$ into the formulation of \cite{keskin2016partial}. Let $K$ denote the set of homogeneous vehicles. We assume that energy consumption rate realization would be different depending on arc $(i,j)$. 

To apply an ARO approach, the decision variables should be categorized into two groups, which are \textit{non-adaptive} (or \textit{here-and-now}) and \textit{adaptive} (or \textit{wait-and-see}) variables. The three-index EVRPTWPR has five types of decision variables: $x_{ijk},\tau_{ik},u_{ik},y_{ik}$, and $Y_{ik}$. Instead of using $Y_{ik}$, we introduce a new decision variable $p_{ik}$, which denotes the energy recharged amount of vehicle $k$ at the station $i$. Among these variables, $x_{ijk}$ and $u_{ik}$ are clearly non-adaptive variables. The route for each vehicle should be decided before the trip, and energy consumption rate uncertainty is not considered. Variables $p_{ik}$ are adaptive variables as they are impacted by the uncertainty of the energy consumption rate. As we explained in Section \ref{ROEVRPTWPR}, the other two decision variables $\tau_{ik}$ and $y_{ik}$ are determined via the energy consumption realization. Therefore, they are adaptive variables.

Instead of minimizing the total distance traveled, we propose to minimize the total energy consumption as an objective function. Since we consider the energy consumption rate as the uncertain parameter, total travel distance may not be the same as minimizing the consumed energy. We note that minimizing the worst energy recharging amount is equivalent to minimizing the worst energy consumption value. Based on this observation, we formulate the AR-EVRPTWPR as follows:
\begin{subequations}
\begin{align}
    \min ~~~& \max_{\boldsymbol{h} \in \mathcal{U}} \sum_{k \in K} \sum_{i \in F_0^{'}} p_{ik} (\boldsymbol{h})
    \label{aroobj}\\
    \textrm{s.t.} ~~~& \sum_{k \in K} \sum_{j \in V_{N+1}^{'}\setminus\{i\}} x_{ijk} = 1 && \hspace{-3cm} \forall i \in V \label{aroconb}\\
    & \sum_{j \in V_{N+1}^{'}, i \neq j} x_{ijk} \le 1 && \hspace{-3cm} \forall i \in F^{'}, k \in K \label{aroconc}\\
    & \sum_{i \in V_0^{'}\setminus\{j\}} x_{ijk} - \sum_{i \in V_{N+1}^{'}\setminus\{j\}} x_{jik} = 0 && \hspace{-3cm} \forall j \in V^{'}, k \in K \label{arocond}\\
    & \sum_{j \in V_{N+1}^{'}} x_{0,j,k} \le 1, && \hspace{-3cm} \forall k \in K \label{arocone}\\
    \notag & \tau_{ik}(\boldsymbol{h}) + (t_{ij} + s_i)x_{ijk} - l_0(1- x_{ijk}) \le \tau_{jk}(\boldsymbol{h}) && \\
    & && \hspace{-3cm} \forall \boldsymbol{h} \in \mathcal{U}, i \in V_0, j \in V_{N+1}^{'}, i \neq j, k \in K \label{aroconf}\\
    \notag & \tau_{ik}(\boldsymbol{h}) + t_{ij}x_{ijk} + g \cdot p_{ik}(\boldsymbol{h}) - (l_0 + gQ)(1 - x_{ijk}) \le \tau_{jk}(\boldsymbol{h}) \\
    & && \hspace{-3cm} \forall \boldsymbol{h} \in \mathcal{U}, i \in F^{'}, j \in V_{N+1}^{'}, i \neq j, k \in K \label{arocong}\\
    & e_j \le \tau_{jk}(\boldsymbol{h}) \le l_j && \hspace{-3cm} \forall \boldsymbol{h} \in \mathcal{U}, j \in V_{0, N+1}^{'}, k \in K \label{aroconh}\\
    & 0 \le u_{jk} \le u_{ik} - q_i x_{ijk} + C(1 - x_{ijk}) && \hspace{-3cm} \forall i \in V_{0}^{'}, j \in V_{N+1}^{'}, i \neq j, k \in K \label{aroconi}\\
    & 0 \le u_0 \le C && \label{aroconj}\\
    \notag & 0 \le y_{jk}(\boldsymbol{h}) \le y_{ik}(\boldsymbol{h}) - (h_{ij} \cdot d_{ij})x_{ijk} + Q(1 - x_{ijk}) && \\
    & && \hspace{-3cm} \forall \boldsymbol{h} \in \mathcal{U}, i \in V, j \in V_{N+1}^{'}, i \neq j, k \in K \label{aroconk}\\
    \notag & 0 \le y_{jk}(\boldsymbol{h}) \le y_{ik}(\boldsymbol{h}) + p_{ik}(\boldsymbol{h}) - (h_{ij} \cdot d_{ij})x_{ijk} + Q(1 - x_{ijk}) \\
    & && \hspace{-3cm} \forall \boldsymbol{h} \in \mathcal{U}, i \in F_0^{'}, j \in V_{N+1}^{'}, i \neq j, k \in K \label{aroconl}\\
    & 0 \le y_{ik}(\boldsymbol{h}) + p_{ik}(\boldsymbol{h}) \le Q && \hspace{-3cm} \forall \boldsymbol{h} \in \mathcal{U}, i \in F_0^{'}, k \in K \label{aroconm}\\
    & 0 \le p_{ik} (\boldsymbol{h}) && \hspace{-3cm} \forall \boldsymbol{h} \in \mathcal{U}, i \in V_{0, N+1}^{'}, k \in K \label{aroconn}\\
    & x_{ijk} \in \{0,1\} && \hspace{-3cm} \forall i \in V_0^{'}, j \in V_{N+1}^{'}, i \neq j, k \in K, \label{arocono}
\end{align}
\label{aromodel}
\end{subequations}
where $\tau(\boldsymbol{h}), y(\boldsymbol{h}), p(\boldsymbol{h})$ denote the functional form of the adaptive decision variables for energy consumption rate scenario $\boldsymbol{h}$. In other words, $\tau, y$, and $p$ are determined for the energy consumption scenarios $\boldsymbol{h}$. The objective (\ref{aroobj}) minimizes the worst recharging amount for all energy consumption rate realization $\boldsymbol{h} \in \mathcal{U}$. Constraints with only non-adaptive decision variables, such as (\ref{aroconb})--(\ref{arocond}) and (\ref{aroconi})--(\ref{aroconj}), are the same as the constraints (\ref{keskinconb})--(\ref{keskincond}) and (\ref{keskinconh})--(\ref{keskinconi}), respectively. We add constraint (\ref{arocone}) to enforce that each vehicle should be dispatched at most once. All of the constraints that contain adaptive decision variables become a robust constraint, such as (\ref{aroconf})--(\ref{aroconh}) and (\ref{aroconk})--(\ref{aroconn}). Constraints (\ref{aroconf})--(\ref{aroconh}) and (\ref{aroconk})--(\ref{aroconn}) enforce arcs' time feasibility, time window, non-negative battery charge level, and bounds of the battery charge level for any energy consumption rate realization $\boldsymbol{h} \in \mathcal{U}$. It is not difficult to see that the presented model is hard to solve using standard methods. In the following section, we introduce a solution framework to solve the proposed model (\ref{aromodel}) and obtain a robust optimal solution.

\section{Solution Methods} \label{solutionmethod}

The AR-EVRPTWPR formulation presented in problem (\ref{aromodel}), depicts a two-stage structure. The first-stage variables are non-adjustable and the second-stage variables are adjustable based on the first stage. In particular, the first stage finds the route and cargo level decision immunized against all possible uncertainty realization $\boldsymbol{h} \in \mathcal{U}$. At the second stage, the decision should be adaptive to any uncertainty realization $\boldsymbol{h} \in \mathcal{U}$, and provide a minimum worst-case cost. To simplify the notation, we present the compact matrix form of the (\ref{aromodel}) as follows. 
\begin{subequations}
\begin{align}
    \min_{\boldsymbol{x}, \boldsymbol{y}(\cdot)} & \max_{\boldsymbol{h} \in \mathcal{U}} ~ \boldsymbol{b}^T \boldsymbol{y}(\boldsymbol{h}) \label{compactobj}\\
    \textrm{s.t.} &~~ \boldsymbol{Ax} \le \boldsymbol{d} \label{compact1stcon}\\
    &~~ \boldsymbol{Fy}(\boldsymbol{h}) \le \boldsymbol{f} \label{compact2ndcon}\\
    &~~ \boldsymbol{G}(\boldsymbol{h})\boldsymbol{x} + \boldsymbol{Hy}(\boldsymbol{h}) \le \boldsymbol{g},\label{compactcomplicated}
\end{align}
\label{compactmodel}
\end{subequations}
where the vector $\boldsymbol{x}$ is the vector of non-adaptive decision variables, including the route and load level decisions. The vector $\boldsymbol{y}$ is the vector of adaptive decision variables, including the arrival time, battery energy level, and battery charge level decisions. The functional form $\boldsymbol{y}(\boldsymbol{h})$ represents that $\boldsymbol{y}$ is fully adaptive to any uncertainty realization $\boldsymbol{h} \in \mathcal{U}$. Constraint (\ref{compact1stcon}) includes route-related constraints (\ref{aroconb})-(\ref{arocone}) and vehicle's cargo level-related constraints (\ref{aroconi}) and (\ref{aroconj}). Constraint (\ref{compact2ndcon}) contains (\ref{aroconh}), (\ref{aroconm}), and (\ref{aroconn}). Constraint (\ref{compactcomplicated}) couples the non-adaptive and adaptive decisions, including (\ref{aroconf}), (\ref{arocong}), (\ref{aroconk}), and (\ref{aroconl}).

We note that the energy consumption rate uncertainty in the AR-EVRPTWPR is column-wise, so robust constraints cannot be dualized to obtain a robust counterpart which is a single-stage. Several solution methods were proposed to solve this type of problem, such as Benders decomposition type cutting plane algorithm (\cite{bertsimas2012adaptive, wang2013two}), linear decision rule approximation (\cite{ben2004adjustable, ben2005retailer, dehghan2017adaptive}), and column-and-constraint generation method (\cite{zeng2013solving}). However, using the linear decision rule, e.g., an affine policy approximation for the adaptive decision variables, is not suitable for the AR-EVRPTWPR because it increases the number of decision variables and the resulting problem easily becomes intractable. The remaining potential solution methods for the AR-EVRPTWPR are the Benders decomposition approach and the column-and-constraint generation method. Both methods are very similar since they decompose the problem into a master and subproblem, and iteratively solve them until the master problem provides the robust optimal solution for all uncertainties. The only difference between them is the type of constraints that are added to the master problem. In the column-and-constraint generation (CCG) method, a copy of variables and constraints of the original model is added to the master problem. On the other hand, single or multi Benders cuts are added to the master problem in the Benders decomposition approach. For the AR-EVRPTWPR, the master problem of the CCG method can be formulated as an EVRPTWPR with multiple energy consumption rate scenarios. It gives us an advantage of evaluating time window and battery capacity violations using existing methods such as the ones presented by \cite{schiffer2018adaptive, cortes2019electric}. Thus, we here adopt a similar approach to the CCG introduced in \cite{zeng2013solving} to solve the AR-EVRPTWPR.

\subsection{Column-and-constraint generation framework} \label{solutionmethod_ccg}

We first introduce the master problem and subproblem of the CCG method. Formulation \eqref{compactmodel} is equivalent to the following model:
\begin{subequations}
\begin{align}
    \min_{\boldsymbol{x}, \eta} ~~& \eta \\
    \textrm{s.t.} ~~& \boldsymbol{Ax} \le \boldsymbol{d} \\
    & \eta \ge \max_{\boldsymbol{h} \in \mathcal{U}} \min_{\boldsymbol{y} \in \boldsymbol{\Omega}(\boldsymbol{x}, \boldsymbol{h})} \boldsymbol{b}^T \boldsymbol{y},
\end{align}
\label{compact2stage}
\end{subequations}
where $\boldsymbol{\Omega}(\boldsymbol{x}, \boldsymbol{h}) = \{\boldsymbol{y} : \boldsymbol{Fy} \le \boldsymbol{f}, \boldsymbol{G}(\boldsymbol{h})\boldsymbol{x} + \boldsymbol{Hy} \le \boldsymbol{g}\}$ is the set of feasible arrival time, battery level, and battery charge level for the fixed route and cargo level decision $\boldsymbol{x}$ and an energy consumption $\boldsymbol{h}$. If the uncertainty set $\mathcal{U}$ has a finite number of energy consumption rate scenarios, then problem (\ref{compact2stage}) can be written as follows:
\vspace{-0.3cm}\begin{subequations}
\begin{align}
    \min_{\boldsymbol{x}, \eta, \boldsymbol{y}^s} ~~& \eta \\
    \textrm{s.t.} ~~& \boldsymbol{Ax} \le \boldsymbol{d} \\
    & \eta \ge \boldsymbol{b}^T \boldsymbol{y}^s \qquad \forall s \in I(\mathcal{U}) \\
    & \boldsymbol{y}^s \in \boldsymbol{\Omega}(\boldsymbol{x}, \boldsymbol{h}^s) \qquad \forall s \in I(\mathcal{U}),
\end{align}
\label{compactccgmodel}
\end{subequations}
where $I(\mathcal{U})$ is the index set of scenarios in the set $\mathcal{U}$, $\boldsymbol{y}^s$ is the vector of decision variable $\boldsymbol{y}$, and $\boldsymbol{h}^s$ is energy consumption scenario $\boldsymbol{h}$ for the $s^{\text{th}}$ scenario in the set $\mathcal{U}$.
%and $\boldsymbol{y}^s$ and $\boldsymbol{h}^s$ are the vector of decision variable $\boldsymbol{y}$ the energy consumption scenario $\boldsymbol{h}$ for the $s^{\text{th}}$ scenario in the set $\mathcal{U}$, respectively.
We note that variable $\boldsymbol{y}$ and the associated constraint set are repeated as much as the number of scenarios in $\mathcal{U}$. The number of constraints and variables in (\ref{compactccgmodel}) %size of problem (\ref{compactccgmodel}) 
is almost $|I(\mathcal{U})|$ times bigger than the deterministic EVRPTWPR, and is thus intractable. To obtain a tractable way to solve this problem, we adopt an iterative approach. CCG starts from the relaxed master problem (RMP), which has a subset of the set $\mathcal{U}$. This RMP can be defined as the EVRPTWPR with multiple scenarios of the energy consumption rate. Let $\bar{\mathcal{U}}$ be the subset of the set $\mathcal{U}$, and $(\boldsymbol{x}^*, \eta^*)$ be the optimal solution of the current RMP. The details of the solution method for solving the RMP are described in Section \ref{solutionmethod_rmp} and \ref{proposed_solutionmethod}. 

The next step identifies the worst-case scenario $\tilde{\boldsymbol{h}}$ defining the largest objective value of the subproblem with a fixed $\boldsymbol{x}^*$. The subproblem is stated as follows:
\begin{subequations}
\begin{align}
    \max_{\boldsymbol{h} \in \mathcal{U}} &\min_{\boldsymbol{y}} ~~\boldsymbol{b}^T \boldsymbol{y} \label{compactsubobj}\\
    \text{s.t.} ~& \boldsymbol{Fy} \le \boldsymbol{f} \label{compactsubcon1}\\
    & \boldsymbol{Hy} \le \boldsymbol{g} - \boldsymbol{G}(\boldsymbol{h})\boldsymbol{x}^*. \label{compactsubcon2}
\end{align}
\label{compactsub}
\end{subequations}
Subproblem (\ref{compactsub}) detects the worst-case scenario $\tilde{\boldsymbol{h}}$ that defines the largest amount of consumed energy for a fixed $\boldsymbol{x}^*$. There could be a scenario that makes the current solution $\boldsymbol{x}^*$ infeasible. To guarantee the feasibility of the subproblem (\ref{compactsub}), slack variables could be added to constraints (\ref{compactsubcon1}) and (\ref{compactsubcon2}). Then, the corresponding penalty terms of slack variables should be added to the objective (\ref{compactsubobj}). The details of a model including the slack variables and the proposed solution approach are described in Section \ref{solutionmehtod_sub}. If optimal value of the subproblem (\ref{compactsub}) is bigger than $\eta^*$, then the detected scenario $\tilde{\boldsymbol{h}}$ is added to the set $\bar{\mathcal{U}}$, and new variables and constraints for $\tilde{\boldsymbol{h}}$ are generated and added to the RMP. The size of the RMP grows every iteration. The overall CCG framework for solving AR-EVRPTWPR is summarized in Algorithm \ref{framework_ccg}. 

\begin{algorithm}
\caption{Column-and-Constraint Generation Framework}
\begin{algorithmic}[1]
\State $\bar{\mathcal{U}} \gets \{\bar{\boldsymbol{h}}\}$
% \State $i \gets 1$
\State $\boldsymbol{x} \gets $ \Call{InitialSolution}{$\bar{\mathcal{U}}$}
\State \textit{RobustFeasibility} $\gets$ \textbf{False}
\While {--\textit{RobustFeasibility}}
    \State $\boldsymbol{x}, \eta \gets $ \Call{RelaxedMasterProblem}{$\boldsymbol{x}, \bar{\mathcal{U}}$} %\Comment{Solve (\ref{compactccgmodel})}
    \State $\boldsymbol{\tilde{h}}, \zeta \gets $ \Call{Subproblem}{$\boldsymbol{x}$} %\Comment{Solve (\ref{compactsub})}
    \If{$\zeta - \eta \ge \delta$} 
        \State $\bar{\mathcal{U}} \gets \bar{\mathcal{U}} \cup \{\boldsymbol{\tilde{h}}\}$
        % \State $i \gets i + 1$
    \Else
        \State \textit{RobustFeasibility} $\gets$ \textbf{True}
    \EndIf
\EndWhile
\end{algorithmic}
\label{framework_ccg}
\end{algorithm}

\subsection{Solution Methods for the Master Problem} \label{solutionmethod_rmp}

As we have seen in the previous section, the RMP of the AR-EVRPTWPR can be formulated as an EVRPTWPR model with multi energy consumption rate scenarios. We refer to this problem as multi-scenario EVRPTWPR. Since the route and cargo level decisions are non-adjustable variables, they should satisfy all scenarios in the given uncertainty set $\bar{\mathcal{U}}$. However, the arrival time, battery energy level, and battery charge level decisions can vary to deal with the different scenarios. We can easily infer that the multi-scenario EVRPTWPR is more challenging to solve than the single-scenario one. To deal with this computational challenge, we adopt one of the solution methods for solving the single-scenario EVRPTWPR, which is a hybrid variable neighborhood search and tabu search (VNS-TS) metaheuristic with an annealing mechanism similar to the one used in \cite{schneider2014electric}. 
{Such a solution method will be presented in detail below, and it has the advantage to combine diversification and intensification features, which are key ingredients in the design of metaheuristics (\cite{Gendreau2010}). The diversification operators are the VNS shaking step (which can generate solutions structurally distant from the incumbent one) and the annealing mechanism (which can accept to move the search to a unimproved neighbor solution). The intensification operator relies on the tabu search metaheuristic that is applied to try to improve the solution provided by the VNS shaking step. The diversification features are able to explore new regions of the solution space, whereas the intensification operator is able to deeply investigate promising regions of the solution space. The combination of intensification/diversification features has been employed successfully in various VRP problems (e.g., \cite{Coindreau2019}).} 
%solution framework presented in \cite{schneider2014electric}. Authors in \cite{schneider2014electric} introduced a hybrid variable neighborhood search and tabu search (VNS-TS) meta-heuristic with an annealing mechanism.

In the following sections, we describe the solution method in details, where all constants and variables in the corresponding penalty function are represented in the time unit. Thus, $h^s_{ij} \gets g \times h^s_{ij} \times d_{ij}$ represents the time needed to recharge the consumed energy on an arc $(i,j)$ for scenario $\boldsymbol{h}^s$, and $H = g \times Q$ denotes the time needed to recharge the vehicle fully.

\subsubsection{Preprocessing, Generalized Cost Function, and Initial Solution}
For computational efficiency, we exclude infeasible arcs as done by \cite{schneider2014electric} and \cite{schiffer2018adaptive}. The infeasibility of arc $(u,v)$ is determined by the following preprocessing constraints:
\begin{subequations}
\begin{align}
    & u, v \in V \wedge q_u + q_v > C, \label{preprocessinga}\\
    & u \in V_0^{'}, v \in V_{N+1}^{'} \wedge e_u + s_u + t_{uv} > l_v, \label{preprocessingb}\\
    & u \in V_0^{'}, v \in V^{'} \wedge e_u + s_u + t_{uv} + s_v + t_{v, N+1} > l_0, \label{preprocessingc} \\
    & u, v \in V \wedge \forall i \in F_0^{'}, j \in F_{N+1}^{'} : h_{iu} + h_{uv} + h_{vj} > H, \label{preprocessingd}\\
    & u, v \in V_{0, N+1}^{'} \wedge h_{uv} > H. \label{preprocessinge}
\end{align}
\label{preprocessing}
\end{subequations}
Constraints (\ref{preprocessinga})-(\ref{preprocessingc}) identify infeasibility based on cargo capacity and time window violations. Constraints (\ref{preprocessingd}) and (\ref{preprocessinge}) identify infeasibility due to the battery capacity as discussed in \cite{schneider2014electric} and \cite{schiffer2018adaptive}, respectively. All infeasible arcs satisfying any of constraints (\ref{preprocessing}) are removed from the arc set $A$. Additionally, at every CCG iteration, we remove infeasible arcs satisfying (\ref{preprocessingd}) and (\ref{preprocessinge}) when a new scenario is added to the set $\bar{\mathcal{U}}$.

As is commonly done, we allow infeasible solutions during the search. The approach of \cite{schneider2014electric} considers the sum of the total traveled distance, total cargo capacity violation, time window violation, battery capacity violation, and the diversification penalty. \cite{schiffer2018adaptive} consider the total recharging time, time window violation, and battery capacity violation. \cite{cortes2019electric} evaluate the cost function as the sum of the time spent at recharging stations, time window violation, and battery capacity violations. In our approach, we define the generalized cost function as follows:
\begin{align}
    f(\mathcal{R}) = \max_{s \in I(\bar{\mathcal{U}})} \bigg\{ \sum_{k \in K} \sum_{i \in F^{'}} p^s_{ik} + \alpha\Big(TW^s(\mathcal{R}) + FL^s(\mathcal{R}) \Big) \bigg\} + \beta FR(\mathcal{R}), \label{general_cost_function}
\end{align}
where $p^s_{ik}, \: TW^s(\mathcal{R})$, and $FL^s(\mathcal{R})$ represent the energy recharged at station $i$ of vehicle $k$, time window violation of routes $\mathcal{R}$, and battery capacity violation of routes $\mathcal{R}$ for the $s^{th}$ scenario in set $\bar{\mathcal{U}}$, respectively. Function $FR(\mathcal{R})$ is the cargo capacity violation of the routes $\mathcal{R}$. All penalty terms are scaled by the weighting factors $\alpha$ and $\beta$. These factors are initialized by $(\alpha_0, \beta_0)$, are dynamically updated during the search, and limited between a given lower bounds $(\alpha_{\min}, \beta_{\min})$ and upper bounds $(\alpha_{\max}, \beta_{\max})$. We note that route $\mathcal{R}$ is said to be \textit{robust feasible} for $\bar{\mathcal{U}}$ if $TW^s(\mathcal{R}) = $ and $FL^s(\mathcal{R}) = $ for all $s \in I(\bar{\mathcal{U}})$ and $FR(\mathcal{R}) = 0$.

To calculate the violations of the time window and battery capacity, we use the corridor-based penalty approach introduced in (\cite{schiffer2018adaptive, cortes2019electric}). The time window and battery capacity violation are calculated as follows:
\begin{subequations}
\begin{align}
    TW^s(\mathcal{R}) & = \sum_{\rho \in \mathcal{R}} \bigg\{ \sum_{v \in \rho} \max \big\{ \min \big\{ a_{s,v}^{\min}, a_{s, v}^{\max} \big\} - l_v, 0 \big\} \bigg\}, && \hspace{-2cm} \forall s \in I(\bar{\mathcal{U}}) \\
    FL^s(\mathcal{R}) & = \sum_{\rho \in \mathcal{R}} \bigg\{ \sum_{v \in \rho} \max \big\{ a_{s, v}^{\min} - a_{s, v}^{\max}, 0 \big\} \bigg\}, && \hspace{-2cm} \forall s \in I(\bar{\mathcal{U}})
\end{align}
\end{subequations}
where $\rho$ denotes the route of each vehicle, $a_{s,v}^{\min}$ is the earliest allowed service start time at a vertex $v$ for scenario $\boldsymbol{h}^s \in \bar{\mathcal{U}}$, and  $a_{s,v}^{\max}$ is the earliest service start time at a vertex $v$ if as much energy as possible is recharged at recharging stations visited before reaching $v$ for scenario $\boldsymbol{h}^s \in \bar{\mathcal{U}}$. Readers can refer to \cite{schiffer2018adaptive} and \cite{ cortes2019electric} for additional details. %in addition to the formulas to calculate violations when a route is updated by insertion or merge. 

The cargo capacity penalty of a route $\rho$ is calculated as follows:
\begin{align}
    FR(\mathcal{R}) = \sum_{\rho \in \mathcal{R}} \max \bigg\{ \sum_{v \in \rho} q_v - C, 0 \bigg\}.
\end{align}
We note that the cargo capacity penalty is calculated in $\mathcal{O}(1)$ for any case where the routes $\mathcal{R}$ is updated (\cite{schneider2014electric, schiffer2018adaptive, cortes2019electric}).

We construct an initial solution using the approach in \cite{lin2021electric}. We first choose an arbitrary point and calculate the angle between the depot, this random point, and the customer. All customers are sorted in ascending order of the angle. Then, the customer is iteratively inserted into an active route at the position causing the smallest energy consumption increment. The route becomes inactive when any customer insertion causes the time window or battery capacity violation, and the other route is activated to insert remaining customers.  The number of available employed vehicles is given as a parameter. If there remain unrouted customers and no more routes can be activated, the remaining customers are inserted into the last route disregarding any time window and batter capacity violation.

\subsubsection{Variable neighborhood search and tabu search components}

The variable neighborhood search (VNS) was proposed by \cite{mladenovic1997variable}, and several other research followed and proposed some variations {(e.g., \cite{Bierlaire2010, schneider2014electric, cortes2019electric, Thevenin2019, lin2021electric}).} Such work showed that VNS could find near optimal solutions, sometimes optimal solutions, for the EVRPTWPR, {and more generally for highly complex optimization problems}. We here introduce VNS-TS proposed by \cite{schneider2014electric} to solve our multi-scenario EVRPTWPR. Given the current solution $\mathcal{R}$, a neighboring solution $\mathcal{R}^{'}$ is generated by a cyclic-exchange operator. It selects a random number $n_r$ of routes from $\mathcal{R}$. An arbitrary length $n_v$ of consecutive vertices is chosen to form an exchange block. These blocks are reversed and transferred between selected routes. The neighboring solution $\mathcal{R}^{'}$ defined by the exchange of blocks could be feasible or not.

{In line with other work on VNS (e.g., \cite{Thevenin2017})}, a tabu search (TS) is used for a local search strategy to improve the neighboring solution $\mathcal{R}^{'}$. For each VNS iteration, TS is run for $n_{tabu}$ iterations, and it provides a $\mathcal{R}^{''}$. If the $\mathcal{R}^{''}$ is an improving solution compared to $\mathcal{R}$, then $\mathcal{R}$ is replaced by $\mathcal{R}^{''}$. However, even if the $\mathcal{R}^{''}$ is not an improving solution, it can be accepted depending on the acceptance criteria based on simulated annealing (SA) (\cite{kirkpatrick1983optimization}). Specifically, the non-improving solution $\mathcal{R}^{''}$ is accepted with a probability of $e^{(f(\mathcal{R}) - f(\mathcal{R}^{''})) / T}$, where $T$ is temperature of annealing phase. Temperature $T$ is initialized by $T_0$ to be accepted non-improving solution $\mathcal{R}^{''}$ with the probability of 50\% when $f(\mathcal{R}^{''})$ is $\Delta_{SA}$ worse than $f(\mathcal{R})$. For every VNS iteration, the temperature linearly decreases. We apply the same parameter settings used by \cite{lin2021electric}.

{In line with the VRP neighborhood structures proposed in the literature (\cite{Golden2008})}, six operators are applied to the TS neighboring solution $\mathcal{R}^{'}$, such as 2-opt$^*$, exchange-intra, exchange-inter, relocate-intra, relocate-inter, and stationInRe. Each operator has its own tabu set. The best none-tabu move is chosen. The move is said to be superior if the resulting solution has a lower cost function value or is a robust feasible solution with less employed vehicles. The operators that are applied are described as follows.

\begin{itemize}
    \item 2-opt$^*$: Select two routes $\rho_i$ and $\rho_j$, and split them into two partial routes, i.e., $\rho_i \rightarrow \rho_i^1, \rho_i^2$ and $\rho_j \rightarrow \rho_j^1, \rho_j^2$. Then, construct two new routes such as $\rho_i^{new} \gets \rho_i^1 + \rho_j^2$ and $\rho_j^{new} \gets \rho_j^1 + \rho_i^2$.
    \item Exchange-intra: Select two vertices $i$ and $j$ from a route $\rho = \langle 0,\dots,i^-,i,\dots,j^-, \allowbreak j,\dots,N+1 \rangle$. Then, exchange the positions of them, such that $\rho^{new} = \langle 0,\dots,i^-, \allowbreak j,\dots,j^-,i,\dots,N+1 \rangle$.
    \item Exchange-inter: Select a vertex $i$ from route $\rho_i$ and a vertex $j$ from route $\rho_j$. Exchange the positions of them, such that $\rho_i^{new} = \langle 0,\dots,i^-,j,i^+\dots,N+1\rangle$ and $\rho_j^{new} = \langle 0,\dots,j^-,i,j^+\dots,N+1\rangle$.
    \item Relocate-intra: Select vertices $i$ and $j$ from route $\rho$. Delete $i$ from $\rho$ and reinsert it after $j$, such that $\rho^{new} = \langle v_0,\dots,i^-,i^+,\dots,j,i,j^+,\dots,v_{N+1} \rangle$.
    \item Relocate-inter: Select a vertex $i$ from route $\rho_i$ and a vertex $j$ from route $\rho_j$. Delete $i$ from $\rho_i$ and reinsert it after $j$ of $\rho_j$, such that $\rho_i^{new} = \langle 0,\dots,i^-,i^+,\dots,N+1 \rangle$ and $\rho_j^{new} = \langle 0,\dots,j,i,j^+,\dots,N+1 \rangle$.
    \item stationInRe: Select an unassigned station vertex and insert it into the route or remove a station vertex in the route.
\end{itemize}
The $i^-$ and $i^+$ represent the preceding and succeeding vertices of the vertex $i$.

Each operator considers all possible vertices and positions. Then, the best move is chosen, and the removed edges by this best move are stored in the corresponding operator's tabu set. Each operator cannot consider the reinsertion of an edge if the edge is in its tabu set. Every edge in the tabu set has its tabu tenure, and the number of tabu tenure is randomly chosen from a range of $\Big[n_{tabu}^{\min}, n_{tabu}^{\max}\Big]$, where $n_{tabu}^{\min}$ and $n_{tabu}^{\max}$ are parameters that will be tuned in the computational results.

\subsection{Solution Method for the Subproblem} \label{solutionmehtod_sub}
Recall the subproblem (\ref{compactsub}) with slack variables and penalty terms presented in Section \ref{solutionmethod_ccg}.
\vspace{-2em}
\begin{subequations}
\begin{align}
    \max_{\boldsymbol{h} \in \mathcal{U}} &\min_{\boldsymbol{y}} ~~\boldsymbol{b}^T \boldsymbol{y} + \boldsymbol{\kappa}^T \boldsymbol{w} \label{compactsubobj_alt}\\
    \text{s.t.} ~& \boldsymbol{Fy} \le \boldsymbol{f} + \boldsymbol{w}_1 \label{compactsubcon1_alt}\\
    & \boldsymbol{Hy} \le \boldsymbol{g} - \boldsymbol{G}(\boldsymbol{h})\boldsymbol{x}^* + \boldsymbol{w}_2, \label{compactsubcon2_alt}
\end{align}
\label{compactsub_alt}
\end{subequations}
where $\boldsymbol{w} = (\boldsymbol{w}_1, \boldsymbol{w}_2)$ and $\boldsymbol{\kappa} = (\boldsymbol{\kappa}_1, \boldsymbol{\kappa}_2)$ the vector of slack variables and penalty coefficients, respectively. As mentioned in Section \ref{solutionmethod_ccg}, slack variables and penalty terms are added to guarantee the feasibility of the subproblem (\ref{compactsub}). Since the subproblem (\ref{compactsub_alt}) has a bilevel structure, it cannot be solved directly. We note that the inner minimization problem is a linear programming problem, and it can be dualized as a maximization problem. Then, the resulting problem can be represented as follows:
\begin{subequations}
\begin{align}
    \max_{\boldsymbol{h}, \boldsymbol{\lambda}, \boldsymbol{\varphi}} ~~~& \boldsymbol{\varphi}^T \big(\boldsymbol{G}(\boldsymbol{h})\boldsymbol{x}^* - \boldsymbol{g} \big) - \boldsymbol{\lambda}^T \boldsymbol{f}  \\
    \textrm{s.t.} ~~& - \boldsymbol{\lambda}^T \boldsymbol{f} - \boldsymbol{\varphi}^T \boldsymbol{H} \le \boldsymbol{b}^T \label{compactinnerdualcon1}\\
    & ~~~ \boldsymbol{0} \le \boldsymbol{\lambda} \le \boldsymbol{\kappa}_1 \\
    & ~~~ \boldsymbol{0} \le \boldsymbol{\varphi} \le \boldsymbol{\kappa}_2 \\
    & ~~~ \boldsymbol{h} \in \mathcal{U}, \label{compactinnerdualcon2}
\end{align}
\label{compactinnerdual}
\end{subequations}
where $\boldsymbol{\lambda}$ and $\boldsymbol{\varphi}$ are multipliers of the constraints (\ref{compactsubcon1}) and (\ref{compactsubcon2}), respectively, and vectors $\boldsymbol{\kappa}_1$ and $\boldsymbol{\kappa}_2$ are the subvectors of $\boldsymbol{\kappa}$ for $\boldsymbol{w}_1$ and $\boldsymbol{w}_2$, respectively. Note that the dualized subproblem (\ref{compactinnerdual}) is a bilinear optimization problem with a bilinear term $\boldsymbol{\varphi}^T \boldsymbol{G}(\boldsymbol{h})\boldsymbol{x}^*$ in the objective function. Since the bilinear terms only exist in the objective function and its feasible region is separable for $(\boldsymbol{\lambda}, \boldsymbol{\varphi})$ and $\boldsymbol{h}$, problem (\ref{compactinnerdual}) is a separable bilinear program. The optimal solution of problem (\ref{compactinnerdual}) consists of the extreme point of the polyhedron $\{(\boldsymbol{\lambda}, \boldsymbol{\varphi}) : - \boldsymbol{\lambda}^T \boldsymbol{f} - \boldsymbol{\varphi}^T \boldsymbol{H} \le \boldsymbol{b}^T, \boldsymbol{0} \le \boldsymbol{\lambda} \le \boldsymbol{\kappa}_1, \boldsymbol{0} \le \boldsymbol{\varphi} \le \boldsymbol{\kappa}_2 \}$ and the extreme point of the uncertainty set $\mathcal{U}$. Therefore, we can solve the subproblem (\ref{compactinnerdual}) using the alternating direction (AD) algorithm proposed by \cite{konno1976cutting}. The AD algorithm always converges to a Karush–Kuhn–Tucker (KKT) point of problem (\ref{compactinnerdual}) (\cite{lorca2014adaptive}). The AD algorithm iteratively solves two problems; one is (\ref{compactinnerdual}) with a fixed $\boldsymbol{h}$ and the other is the same problem (\ref{compactinnerdual}) with a fixed $(\boldsymbol{\lambda}, \boldsymbol{\varphi})$. Details are summarized in Algorithm \ref{AD_sub}.
\begin{algorithm}
\caption{Alternating Direction (AD) Algorithm}
\begin{algorithmic}[1]
\Procedure{Subproblem}{$\boldsymbol{x}^*$}
\State $\tilde{\boldsymbol{h}} \gets \bar{\boldsymbol{h}}$
\State $LB \gets 0$, $UB \gets \infty$
\While {$UB - LB < \varepsilon$}
    \State $LB \gets \max_{(\boldsymbol{\lambda}, \boldsymbol{\varphi}) \in \Pi} \boldsymbol{\varphi}^T \big(\boldsymbol{G}(\tilde{\boldsymbol{h}})\boldsymbol{x}^* - \boldsymbol{g} \big) - \boldsymbol{\lambda}^T \boldsymbol{f}$ and let $(\boldsymbol{\lambda}^*, \boldsymbol{\varphi}^*)$ be its optimal solution 
    \Statex \Comment{$\Pi := \{(\boldsymbol{\lambda}, \boldsymbol{\varphi}) : \text{(\ref{compactinnerdualcon1})--(\ref{compactinnerdualcon2})} \}$}
    \State $UB \gets \max_{\boldsymbol{h} \in \mathcal{U}} \boldsymbol{\varphi}^{*T} \big(\boldsymbol{G}(\boldsymbol{h})\boldsymbol{x}^* - \boldsymbol{g} \big) - \boldsymbol{\lambda}^{*T} \boldsymbol{f}$ and let $\boldsymbol{h}^*$ be its optimal solution
    \State $\tilde{\boldsymbol{h}} \gets \boldsymbol{h}^*$
\EndWhile
\State $\zeta \gets UB$
\EndProcedure
\end{algorithmic}
\label{AD_sub}
\end{algorithm}

\subsection{Proposed Solution Framework} \label{proposed_solutionmethod}

In the previous sections, we introduce the CCG method, the VNS-TS algorithm, and the AD algorithm. The general CCG method iteratively solves the master problem to optimality and detects a scenario defining the worst violation. In other words, multi-scenario EVRPTWPR has to be solved several times, and the size of the problem increases as the CCG iteration number increases. Therefore, we can easily infer that our multi-scenario EVRPTWPR suffers computational difficulties compared to the single-scenario EVRPTWPR. To avoid an expected computational intractability, we propose a CCG framework motivated by the Benders based branch-and-bound (BBC) framework as presented by \cite{naoum2013interior}. The BBC algorithm finds Benders cuts from the incumbent solution, which is not optimal for the RMP. More precisely, if a node in the branch-and-bound tree of the master problem specifies a feasible solution, then the subproblem is solved, and Benders cuts are added to the master problem. Previous work showed that this approach can find the optimal solution efficiently in a decomposition framework (\cite{adulyasak2015benders, lorca2016multistagedynamic}).

The proposed solution method starts from an uncertainty set with only a nominal energy consumption rate scenario $\bar{\mathcal{U}} = \{\bar{\boldsymbol{h}}\}$. The VNS-TS algorithm is used to find a robust feasible route for $\bar{\mathcal{U}}$. If any robust feasible route is found, then the subproblem is solved to detect a scenario defining the worst violation. The detected scenario is added to the uncertainty set $\bar{\mathcal{U}}$. Once any scenario is added, the current robust feasible route can be infeasible or have higher cost than before. The VNS-TS keeps going to find a robust feasible route for the current uncertainty set or improve the current route. The stopping criteria of the VNS-TS algorithm is that route $\mathcal{R}$ should be robust feasible for the uncertainty set $\bar{\mathcal{U}}$ and that the gap between the value of the generalized cost function $f(\mathcal{R})$ and the objective value of the subproblem should be less than the predefined threshold. We summarize the proposed algorithm in Algorithm \ref{framework_ccg_vns}.

\begin{algorithm}
\caption{Proposed Algorithm Based on CCG framework}

\begin{algorithmic}[1]
\State $\bar{\mathcal{U}} \gets \{\bar{\boldsymbol{h}}\}$
\State $i \gets 0 $
\State $UB \gets \infty$
\State $x \gets $ generateInitialSolution()
\While {\textbf{True}}
    \State $x' \gets $ \Call{cyclic-exchange}{$x$}
    \State $x''\gets $ \Call{TabuSearch}{$x'$}
    \If{$f(x'') \le UB$}
        \State $i \gets 0$
        \State $x \gets  x''$
        \If{$x''$ is robust feasible for $\mathcal{\bar{U}}$}
            \State $\zeta, \boldsymbol{h}^* \gets $ \Call{AlternatingDirection}{$x''$}
            \If{$\zeta \le f(x'')$}
                \State $x''$ is the best robust feasible solution for $\mathcal{U} \rightarrow$ Terminate Algorithm
            \Else
                \State $UB \gets f(x'')$
                \State $\mathcal{\bar{U}} \gets \mathcal{\bar{U}} \cup \{\boldsymbol{h}^*\}$
            \EndIf
        \EndIf
    \ElsIf{\Call{acceptSA}{$x, x''$}}
        \State $x \gets  x''$
    \EndIf
    \If{$x''$ is robust feasible for $\mathcal{\bar{U}}$}
        \State $\zeta, \boldsymbol{h}^* \gets $ \Call{AlternatingDirection}{$x''$}
        \If{$\zeta - f(x'') < \delta$}
            \State $x''$ is the best robust feasible solution for $\mathcal{U} \rightarrow$ Terminate Algorithm
        \EndIf
    \Else
        \State $i \gets i + 1$
        \If{$i > n_{feas}$}
            \State No robust feasible solution $\rightarrow$ Terminate Algorithm
        \EndIf
    \EndIf
\EndWhile

\end{algorithmic}
\label{framework_ccg_vns}
\end{algorithm}

We note that detecting the worst violation scenario whenever a robust feasible route is found is similar to that when the BBC framework finds Benders cuts from the incumbent solution. Instead of using the VNS-TS algorithm to find the best solution for the multi-scenario EVRPTWPR, we use the VNS approach to play the role of the branch-and-bound tree in the BBC framework. Also, any robust feasible route for the uncertainty set $\bar{\mathcal{U}}$ in the VNS-TS procedure can be interpreted as the incumbent solution in the BBC framework. 

Since the uncertainty set $\bar{\mathcal{U}}$ is a subset of the set $\mathcal{U}$, any robust feasible solution of the RMP with $\bar{\mathcal{U}}$ provides the lower bound of the AR-EVRPTWPR. We also note that the subproblem provides the upper bound. Therefore, the stopping criterion is sufficient to guarantee that the route obtained from our proposed solution method is robust feasible for all uncertainty scenarios in the set $\mathcal{U}$. 

\section{Computational experiments}\label{results}

This section presents extensive computational results to illustrate the performance of the proposed AR-EVRPTWPR. We first show the results related to the proposed Algorithm 1 for the AR-EVRPTWPR. Then, Monte-Carlo simulation tests are conducted to evaluate the economical impact and robustness of the route of each model.

The AR-EVRPTWPR is first applied to small-sized instances of \cite{schneider2014electric} to demonstrate the quality of solutions obtained using our algorithm. We compare these results with the results from the commercial solver CPLEX. Next, we generate mid-sized instances to show the computational capability of our algorithm. These mid-sized instances are generated as follows. We randomly pick 20, 30, and 40 customers are from the large-sized instances of \cite{schneider2014electric}. Then, we solve the created mid-sized instances using the VNS-TS {metaheuristic} five times, and the visited stations in the shortest route are chosen. We note that this instance generation is similar to how \cite{schneider2014electric} generated the small-sized instances. %But the only subset of the large-sized instances used in small-sized instances creation is considered, and recharging stations are chosen from the shortest route of five times run to avoid visiting too many stations.

For the uncertainty set, we consider the following cardinality constrained uncertainty set proposed by \cite{bertsimas2004price}:
\begin{align}
    \mathcal{U} := \Bigg\{ \boldsymbol{h} \in \mathbb{R}^{|\mathcal{A}|} : \sum_{a \in \mathcal{A}} \frac{|h_a - \bar{h}_a|}{\hat{h}_a} \le \Gamma, h_a \in [\bar{h}_a - \hat{h}_a, \bar{h}_a + \hat{h}_a], \forall a \in \mathcal{A} \Bigg\} \label{uncertaintyset_budget}
\end{align}
where $\bar{h}_a$ is the expected energy consumption rate of arc $a$, $\hat{h}_a := 0.1\times \bar{h}_a$, and $\Gamma$ is the \enquote{budget of uncertainty} (\cite{bertsimas2012adaptive}). The uncertainty set (\ref{uncertaintyset_budget}) assumes that the energy consumption rate realizes in a certain range. The cardinality constraint limits the total deviation of uncertain energy consumption rates from their expected values. Since the purpose of this paper is to propose an ARO model of EVRPTWPR and show its advantage, we here did not perform an extensive sensitivity analysis for the different Gamma values. The following sections will show that considering uncertainty could provide significant results even if the value of $\Gamma$ is small. %We left various types of uncertainty sets and their sensitivity analysis as future work.

\subsection{Experimental setup}
All models and algorithms are implemented in Python and performed on a server with an Intel Xeon Silver 4114 CPU. CPLEX 20.1.0 is used as the optimization solver.

For the parameter settings, all CPLEX parameters are kept to the default values. The convergence tolerance for the CCG framework is set as $\delta = 0.001$. The penalty cost vector $\boldsymbol{\kappa}$ of the subproblem is set as $10000 \cdot \boldsymbol{1}$, where $\boldsymbol{1}$ is the vector where all elements have a value of 1. Since the objective value of the subproblem can be huge in the case of infeasibility, the stopping criterion for the AD algorithm is set as $(UB - LB) / LB < 0.0001$. For the VNS-TS algorithm, the number of VNS iterations is limited to $n_{feas} = 500$ for all sizes of instances. The penalty parameters $(\alpha, \beta)$ are set to $(\alpha_{0}, \beta_{0}) = (10,10)$, $(\alpha_{min}, \beta_{min}) = (1,1)$ and $(\alpha_{max}, \beta_{max}) = (5000,5000)$. The penalty parameters increase or decrease by 1.2 times every two continuous infeasible and feasible VNS iterations, respectively. For the simulated annealing based acceptance criteria, $\Delta_{SA}$ is set to 10. For the tabu search, the lower and upper bounds of tabu tenure are set as $n_{tabu}^{min} = 15$ and $n_{tabu}^{max} = 30$, respectively.

A set of homogeneous EVs is considered, where each EV has a battery capacity $Q$ equals to 77.75, cargo capacity $C$ of 200, and recharging rate $g$ of 1.0. In addition, the expected energy consumption rate $\bar{h}_a$ is set to 1.0 for all arcs between vertices. We note that all these values are the same as the small-sized instances of \cite{schneider2014electric}. In addition, we allow each vehicle to re-visit the same recharging station only once.

\subsection{Computational results}

The computational performance of the proposed Algorithm \ref{framework_ccg_vns} is reported in Table \ref{table_cplex_ccg}. We compare CCG using CPLEX with Algorithm \ref{framework_ccg_vns} in terms of solution time and quality for small-sized instances of \cite{schneider2014electric}. For CCG with CPLEX, the RMP is solved to optimality using CPLEX for each CCG iteration, and the AD algorithm solves the subproblem as presented in Algorithm \ref{AD_sub}. The computational time is limited to 10,800 CPU seconds for both the CCG with CPLEX and for Algorithm \ref{framework_ccg_vns} and $\Gamma$ is set to 6. We note that there are four small-sized instances where the robust feasible solution does not exist for the given uncertainty set ($\hat{h}_a = 0.1 \bar{h}_a, ~ \forall a \in \mathcal{A}$ and $\Gamma = 6$). There is one in the 5-customers instance (i.e., R105C5) and three in the 10-customers instance (i.e., R103C10, R203C10, and RC102C10).

In Table \ref{table_cplex_ccg}, column \enquote{$Inst.$} represents the instance's name, and the last value in each instance's name represent the number of customers (i.e., 5, 10, and 15). Column \enquote{$n_v$} denotes the number of employed vehicles. Columns \enquote{$worst$} and \enquote{$dist.$} records the worst-case energy consumption and the total distance of routes, respectively. Column \enquote{$t(s)$} reports the total CPU time in seconds. Finally, columns \enquote{$gap_w$} and \enquote{$gap_d$} represent the increment of Algorithm \ref{framework_ccg_vns} compared to the results of the CCG with CPLEX in terms of worst-case energy consumption and the total distance of routes, respectively. Both of them are reported in percentage. The \enquote{10800*} in the $t(s)$ column denotes that the algorithm fails to terminate within the time limit. For these cases, we report the results of the best robust feasible solution found during searching. Finally, \enquote{T} indicates that the algorithm fails to find any robust feasible solution within the time limit.

Results in Table \ref{table_cplex_ccg} show the proposed Algorithm \ref{framework_ccg_vns} is promising to solve the AR-EVRPTWPR. The CCG framework using CPLEX can solve all of the instances with 5-customers. On the other hand, it can find robust feasible solutions for 7 out of the 9 10-customers instances and fails to find any robust feasible solution for most instances with 15 customers. However, our Algorithm \ref{framework_ccg_vns} can find robust feasible solutions for all small-sized instances within a reasonable computational time. Except for only one instance (i.e., C202C10), the proposed Algorithm \ref{framework_ccg_vns} finds robust optimal or near-optimal solutions. Also, it can provide better solutions compared to the robust feasible solutions that the CCG with CPLEX finds within the time limit. For instance, the best robust feasible solution found by the proposed Algorithm \ref{framework_ccg_vns} has 24.65\% shorter travel distance (at RC204C15) than the solution provided by the CCG with CPLEX  and the computational time was only 69.22 seconds.

Table \ref{table_mid_vns} reports the computational results for mid-sized instances. Column \enquote{Deterministic model} represents the number of employed vehicles and travel distances when energy consumption rate uncertainty is not considered. Column \enquote{ARO model} denotes the results of AR-EVRPTWPR consisting of the number of employed vehicles, worst-case energy consumption, travel distance, and computational time. For the uncertainty set, the $\Gamma$ is set to 6 and $\hat{h}_a = 0.1\bar{h}_a, ~ \forall a \in \mathcal{A}$. The computational time limit is set to 10,800 CPU seconds. 

Table \ref{table_mid_vns} illustrates that the proposed Algorithm \ref{framework_ccg_vns} can solve all mid-sized instances within the time limit. The computational time increases as the instance's size increases, and most of the computational time is less than two hours. Results in Table \ref{table_mid_vns} also show that the routes become more conservative than the deterministic ones. In most cases, the travel distances increase. Also, for some instances, the required number of vehicles is higher than the deterministic model. It can be interpreted as follows: to withstand the energy consumption rate uncertainties, the employed vehicles should visit more charging stations or fewer customers than the deterministic cases.

From the results of Tables \ref{table_cplex_ccg} and \ref{table_mid_vns}, we can demonstrate that the proposed Algorithm \ref{framework_ccg_vns} is able to find a robust feasible solution that is near to the robust optimal solution. Furthermore, the computational time is reasonable up to 40-customers cases, i.e., less than two hours. 

\begin{table}
\small
\renewcommand{\arraystretch}{0.8}
\centering 
\caption{Comparing CCG using CPLEX with the proposed approach}
\vspace{3mm}
\begin{tabular}{c c c c c c c c c c c}
\toprule
\rule{0pt}{2.7ex}
{} & \multicolumn{4}{c}{CCG with CPLEX} & \multicolumn{6}{c}{Proposed Algorithm \ref{framework_ccg_vns}} \\
\cmidrule(lr){2-5} \cmidrule(lr){6-11}
\vspace{0.3em}{$Inst.$} & {$n_v$} & {$worst$} & {$dist.$} & {$t(s)$} & {$n_v$} & {$worst$} & {$dist.$} & {$t(s)$} & {$gap_{w}$} & {$gap_d$} \\
\hline
\rule{0pt}{2.7ex}
C101C5	&	2	&	277.47	&	257.75	&	219.89	&	2	&	277.47	&	257.75	&	0.48	&	0	&	0	\\ % result of (2, 0)
C103C5	&	1	&	190.02	&	176.12	&	4.10	&	1	&	190.02	&	176.12	&	0.39	&	0	&	0	\\ % result of (1, 1)
C206C5	&	1	&	269.36	&	250.28	&	2121.17	&	1	&	269.36	&	250.28	&	0.73	&	0	&	0	\\ % result of (1, 1)
C208C5	&	1	&	179.88	&	164.34	&	205.35	&	1	&	179.88	&	164.34	&	0.41	&	0	&	0	\\ % result of (1, 1)
R104C5	&	2	&	148.00	&	137.25	&	171.68	&	2	&	148.00	&	137.25	&	0.34	&	0	&	0	\\ % result of (2, 0)
R202C5	&	1	&	154.66	&	143.29	&	216.54	&	1	&	154.66	&	143.29	&	0.39	&	0	&	0	\\ % result of (1, 0)
R203C5	&	1	&	193.13	&	179.06	&	81.97	&	1	&	193.13	&	179.06	&	0.32	&	0	&	0	\\ % result of (1, 0)
RC105C5	&	2	&	260.58	&	241.30	&	118.29	&	2	&	260.58	&	241.30	&	3.79	&	0	&	0	\\ % result of (2, 0)
RC108C5	&	2	&	274.38	&	253.93	&	13.79	&	2	&	274.38	&	253.93	&	0.40	&	0	&	0	\\ % result of (2, 0)
RC204C5	&	1	&	191.99	&	176.40	&	70.67	&	1	&	191.99	&	176.40	&	0.38	&	0	&	0	\\ % result of (1, 0)
\vspace{0.3em} RC208C5	&	1	&	189.42	&	173.57	&	12.53	&	1	&	189.42	&	173.57	&	1.13	&	0	&	0	\\ % result of (1, 0)
\hline																		\rule{0pt}{2.7ex}
C101C10	&	3	&	415.82	&	393.56	&	9024.16	&	3	&	415.82	&	393.56	&	46.63	&	0	&	0	\\ % result of (3,0)
C104C10	&	2	&	304.25	&	287.84	&	10800*	&	2	&	290.21	&	273.93	&	22.02	&	-4.61	&	-4.83	\\ % result of (2,0)
C202C10	&	2	&	266.88	&	251.95	&	10800*	&	2	&	267.21	&	252.05	&	24.88	&	0.12	&	0.04	\\ % result of (,)
C205C10	&	2	&	245.18	&	231.97	&	10800*	&	2	&	245.18	&	231.97	&	67.99	&	0	&	0	\\ % result of (,)
R102C10	&	4	&	275.96	&	262.93	&	10800*	&	4	&	275.96	&	262.93	&	2.67	&	0	&	0	\\ % result of (,)
R201C10	&	2	&	240.19	&	228.41	&	7226.32	&	2	&	240.19	&	228.41	&	26.79	&	0	&	0	\\ % result of (,)
RC108C10	&	\multicolumn{4}{c}{T}		&	3	&	378.71	&	357.23	&	13.74	&	--	&	--	\\ % result of (,)
RC201C10	&	2	&	385.67	&	363.53	&	10800*	&	2	&	377.27	&	355.89	&	15.99	&	-2.18	&	-2.10	\\ % result of (,)
\vspace{0.3em} RC205C10	&	\multicolumn{4}{c}{T}		&	2	&	391.42	&	368.99	&	1.56	&	--	&	--	\\ % result of (2, 0)
\hline																		\rule{0pt}{2.7ex}
C103C15	&	\multicolumn{4}{c}{T}		&	3	&	396.16	&	380.14	&	98.80	&	--	&	--	\\ % result of (,)
C106C15	&	\multicolumn{4}{c}{T}		&	3	&	378.75	&	361.53	&	73.58	&	--	&	--	\\ % result of (,)
C202C15	&	\multicolumn{4}{c}{T}		&	3	&	394.48	&	377.95	&	217.01	&	--	&	--	\\ % result of (,)
C208C15	&	2	&	322.30	&	306.25	&	10800*	&	2	&	318.04	&	300.55	&	28.02	&	-1.32	&	-1.86	\\ % result of (2, 1)
R102C15	&	\multicolumn{4}{c}{T}		&	5	&	456.66	&	441.13	&	648.25	&	--	&	--	\\ % result of (,)
R105C15	&	\multicolumn{4}{c}{T}		&	4	&	368.63	&	353.93	&	74.28	&	--	&	--	\\ % result of (,)
R202C15	&	\multicolumn{4}{c}{T}		&	3	&	401.33	&	386.28	&	138.33	&	--	&	--	\\ % result of (,)
R209C15	&	2	&	420.43	&	403.51	&	10800*	&	2	&	322.07	&	307.68	&	65.97	&	-23.40	&	-23.75	\\ % result of (2, 1)
RC103C15	&	\multicolumn{4}{c}{T}		&	4	&	416.26	&	397.75	&	114.84	&	--	&	--	\\ % result of (,)
RC108C15	&	\multicolumn{4}{c}{T}		&	3	&	399.11	&	378.35	&	26.81	&	--	&	--	\\ % result of (,)
RC202C15	&	\multicolumn{4}{c}{T}		&	3	&	418.01	&	397.20	&	382.40	&	--	&	--	\\ % result of (,)
RC204C15	&	2	&	432.52	&	412.47	&	10800*	&	2	&	327.65	&	310.81	&	69.22	&	-24.25	&	-24.65	\\ % result of (2, 1)
\bottomrule
\end{tabular}
\label{table_cplex_ccg}
\end{table}

\begin{table}
\small
\renewcommand{\arraystretch}{0.8}
\centering 
\caption{Results for mid-sized instances}
\vspace{3mm}
\begin{tabular}{c c c c c c c}
\toprule
\rule{0pt}{2.7ex}
{} & \multicolumn{2}{c}{Deterministic model} & \multicolumn{4}{c}{ARO model} \\
\cmidrule(lr){2-3} \cmidrule(lr){4-7}
\vspace{0.3em}{$Inst.$} & {$n_v$} & {$dist.$} & {$n_v$} & {$worst$} & {$dist.$} & {$t(s)$}\\
\hline
\rule{0pt}{2.7ex}
C101C20	&	4	&	466.63	&	4	&	483.50	&	466.82	&	284.34	\\
C103C20	&	4	&	381.83	&	4	&	399.85	&	383.25	&	292.53	\\
C206C20	&	2	&	412.79	&	2	&	445.18	&	425.61	&	233.35	\\
C208C20	&	2	&	376.61	&	2	&	491.58	&	472.31	&	205.88	\\
R104C20	&	4	&	337.26	&	5	&	378.64	&	364.36	&	187.87	\\
R105C20	&	6	&	505.07	&	7	&	568.12	&	549.28	&	230.75	\\
R203C20	&	2	&	260.26	&	2	&	273.78	&	262.94	&	164.60	\\
RC105C20	&	8	&	654.70	&	8	&	685.80	&	664.26	&	232.39	\\
RC108C20	&	6	&	520.83	&	6	&	569.25	&	546.25	&	1065.22	\\
\vspace{0.3em} RC204C20	&	2	&	410.26	&	4	&	428.56	&	409.91	&	174.72	\\
\hline
\rule{0pt}{2.7ex}
C101C30	&	5	&	739.64	&	5	&	772.66	&	750.45	&	5739.58	\\
C104C30	&	4	&	459.00	&	4	&	476.02	&	460.69	&	750.46	\\
C202C30	&	3	&	516.28	&	3	&	548.84	&	530.79	&	761.65	\\
C205C30	&	3	&	539.84	&	3	&	560.57	&	540.56	&	463.27	\\
R102C30	&	7	&	604.75	&	8	&	603.36	&	585.64	&	1214.53	\\
R103C30	&	8	&	740.98	&	9	&	790.43	&	774.16	&	4715.54	\\
R201C30	&	3	&	567.26	&	3	&	626.90	&	609.70	&	4953.61	\\
R203C30	&	3	&	404.22	&	3	&	434.10	&	417.03	&	1099.81	\\
RC102C30	&	10	&	904.62	&	11	&	1031.00	&	999.94	&	5916.7	\\
\vspace{0.3em} RC205C30	&	4	&	629.98	&	4	&	663.83	&	644.97	&	947.61	\\
\hline
\rule{0pt}{2.7ex}
C103C40	&	5	&	685.23	&	5	&	729.52	&	712.29	&	2868.03	\\ % duplicate necessary
C202C40	&	4	&	516.01	&	4	&	548.24	&	534.29	&	3935.44	\\ % duplicate necessary
C208C40	&	4	&	524.94	&	4	&	539.22	&	525.31	&	4084.10	\\
R102C40	&	11	&	917.60	&	12	&	968.35	&	947.42	&	2479.40	\\ % duplicate necessary
R105C40	&	10	&	821.08	&	11	&	905.94	&	887.3	&	932.01	\\ % no duplicate required
R202C40	&	4	&	570.40	&	4	&	600.98	&	585.45	&	1179.85	\\ % duplicate necessary
R209C40	&	4	&	499.92	&	4	&	515.84	&	503.04	&	3194.52	\\ % duplicate necessary
RC108C40	&	10	&	874.35	&	11	&	955.82	&	933.09	&	3858.94	\\ % duplicate necessary
RC202C40	&	4	&	834.36	&	5	&	761.95	&	741.65	&	2759.23	\\ % duplicate necessary
RC204C40	&	4	&	598.25	&	4	&	622.11	&	605.76	&	1534.38	\\ % duplicate necessary
\bottomrule
\end{tabular}
\label{table_mid_vns}
\end{table}

\subsection{Economical impact and robustness} \label{section_simulation}

This section studies the economic efficiency and robustness of the ARO model by conducting Monte-Carlo simulation tests for 40-customer instances. We generated a set of 1000 energy consumption rate scenarios using three probability distributions: normal distribution (ND), uniform distribution (UD), and normal distribution considering correlation (NDC). For ND and UD, the energy consumption rate of arc $a$, $h_a$, is sampled from the normal distribution with $\mu = \bar{h}_a, \sigma = 0.1 \bar{h}_a$ and the uniform distribution with the range $[0.9\bar{h}_a, 1.1\bar{h}_a]$, respectively. For NDC, \cite{lee2012robust} assume that customers in congested areas tend to have more travel time than customers in other areas. Similarly, we postulate that arcs located close to one another have similar energy consumption rate uncertainty. Thus, $h_a$ is defined as $\bar{h}_a + \bar{h}_a (z_i + z_j)$, where $z_a$ follows $\mathcal{N}(0, (\frac{0.1}{2\sqrt{2}})^2)$ for all arcs $a \in \mathcal{A}$.

Since deterministic or robust routes can be infeasible for some scenarios, we use the generalized cost function $f(\mathcal{R})$ in (\ref{general_cost_function}) to evaluate the route's costs. In the simulation tests, $\alpha$ and $\beta$ in the generalized cost function $f(\mathcal{R})$ are set to $(10,10)$. For example, if route $\mathcal{R}$ is feasible for scenario $S1$, its cost will be electric energy to complete EVs' travel (i.e., $f(\mathcal{R}) = \sum_{k \in K, i \in F^{'}} p_{ik}$). However, if route $\mathcal{R}$ is infeasible for scenario $S2$, its cost will be the sum of required energy and penalties (i.e., $f(\mathcal{R}) = \sum_{k \in K, i \in F^{'}} p_{ik} + 10 \times ( TW^{S2} (\mathcal{R}) + FL^{S2} (\mathcal{R}) + FR(\mathcal{R}) )$).

Tables \ref{table_simul_normal_correlated}, \ref{table_simul_normal}, and \ref{table_simul_uniform} report the simulation results for NDC, ND, and UD, respectively. Column \enquote{$\Gamma$} represents the considered value of $\Gamma$ in a model. A model with $\Gamma = 0$ denotes the deterministic model. Columns \enquote{Avg.} and \enquote{Std.} represent the average route's costs and the standard deviation of the route's costs. Column \enquote{Max} and \enquote{Min} represent the largest and smallest route's cost among 1000 scenarios. To show how the budget of uncertainty affects the robustness, we report the results of two ARO models with $\Gamma = 6$ and $\Gamma = 12$ and compare them with the deterministic model. In addition, there are some cases where the ARO model employs more EVs. For those instances and to have a fair comparison, we also report the results of a deterministic model that uses the same number of EVs.

The results in Tables \ref{table_simul_normal_correlated}, \ref{table_simul_normal}, and \ref{table_simul_uniform} show the following. First, a route from the ARO model has a higher capability to complete their service. For all probability distributions, ARO models' results indicate higher feasibility than the deterministic model. For instance, in the case of C103C40 and uniform distribution, the deterministic model is only feasible in 35.4\%. Still, the ARO model with Gamma=12 is feasible for most scenarios except two. Also, even if the ARO model's feasibility is low, it is higher than the deterministic model from 2.2 (as seen for R102C40 instance in Table \ref{table_simul_normal_correlated}) to 6.5 times (for the same R102C40 instance in Table \ref{table_simul_normal}). Second, the results of the deterministic model show higher volatility than ARO models. We can see that the ARO models show a lower std for all cases, and its magnitude is in the range from 1.1 (R209C40 in Table \ref{table_simul_normal_correlated}) to 13.2 times (RC202C40 in Table \ref{table_simul_uniform}). The high value of the std can be interpreted as the insecurity of the route. In other words, there could be more risk of failure to complete travelling or paying more charging costs than expected. This pattern can also be observed by comparing the different Gamma values. The ARO models with $\Gamma=12$ usually have lower std values than the ARO models with $\Gamma=6$. Lastly, more employed EVs cannot guarantee robustness. For instances R102C40, R105C40, RC108C40, and RC202C40, ARO models employ one additional EV to find their robust solution. Except for the instance of RC202C40, the other three instances show that the deterministic model provides similar results whether an additional EV is employed or not. In RC202C40, even if additional EVs can provide a more robust route, the results are still worse than the ARO models. For instance, a route employing 5 EVs of the deterministic model shows 71\% feasibility and 57.38 MWh for std, but a route of ARO model with $\Gamma = 12$ shows 100\% feasibility and 10 MWh for std (see Table \ref{table_simul_normal_correlated}). 

\begin{table}
\small
\renewcommand{\arraystretch}{0.8}
\centering 
\caption{Simulation Results - Normal Distribution Considering Correlation}
\vspace{3mm}
\begin{tabular}{c c c c c c c c}
\toprule
\vspace{0.3em}
{$Inst.$} & {$n_v$} & {$\Gamma$} & {Avg.} & {Std.} & {Max} & {Min} & {Feasibility(\%)} \\
\hline
\rule{0pt}{2.7ex}
C103C40	&   5   &	0	&	739.31	&	76.06	&	1181.19	&	654.88	&	37.8	\\
	    &   5   &	6	&	735.72	&	41.18	&	923.20	&	679.99	&	54.0	\\
\vspace{0.5em}	    &   5   &	12	&	715.40	&	15.11	&	887.35	&	682.68	&	97.2	\\
C202C40	&   4   &	0	&	526.05	&	29.45	&	716.56	&	493.75	&	79.6	\\
    	&   4   &	6	&	534.63	&	8.03	&	590.34	&	510.83	&	99.4	\\
\vspace{0.5em}	    &   4   &	12	&	553.94	&	7.60	&	577.61	&	528.94	&	100	\\
C208C40	&   4   &	0	&	524.56	&	6.93	&	548.57	&	505.15	&	100	\\
    	&   4   &	6	&	524.96	&	7.11	&	548.14	&	503.18	&	100	\\
\vspace{0.5em}    	&   4   &	12	&	530.86	&	7.07	&	555.14	&	510.39	&	100	\\
R102C40	&   11   &	0	&	1050.95	&	184.36	&	2124.14	&	854.53	&	24.2	\\
        &   12   &	0	&	1043.61	&	156.06	&	1901.64	&	867.58	&	27.2	\\
        &   12   &	6	&	986.83	&	80.95	&	1430.77	&	878.82	&	53.6	\\
\vspace{0.5em}        &   12   &	12	&	968.18	&	60.97	&	1341.97	&	879.75	&	59.1	\\
R105C40	&   10   &	0	&	901.28	&	110.71	&	1600.49	&	767.45	&	29.0	\\
        &   11   &	0	&	911.00	&	116.12	&	1701.84	&	772.83	&	28.7	\\
        &   11   &	6	&	907.21	&	55.27	&	1327.24	&	827.66	&	68.9	\\
\vspace{0.5em}        &   11   &	12	&	911.10	&	45.08	&	1173.85	&	841.42	&	70.9	\\
R202C40	&   4   &	0	&	606.05	&	55.75	&	903.26	&	547.75	&	49.0	\\
        &   4   &	6	&	587.89	&	13.56	&	703.34	&	562.71	&	93.7	\\
\vspace{0.5em}        &   4   &	12	&	589.04	&	14.36	&	724.95	&	559.99	&	94.6	\\
R209C40	&   4   &	0	&	500.48	&	8.91	&	593.97	&	476.77	&	98.2	\\
    	&   4   &	6	&	503.43	&	8.17	&	585.27	&	480.85	&	98.8	\\
\vspace{0.5em}    	&   4   &	12	&	524.19	&	7.97	&	619.76	&	499.61	&	99.4	\\
RC108C40	&   10   &	0	&	979.76	&	131.23	&	1559.45	&	816.64	&	26.0	\\
            &   11   &	0	&	1025.77	&	158.27	&	1769.91	&	827.88	&	19.8	\\
            &   11   &	6	&	958.80	&	65.32	&	1316.57	&	860.96	&	63.4	\\
\vspace{0.5em}            &   11   &	12	&	1006.53	&	60.99	&	1356.31	&	911.71	&	67.7	\\
RC202C40	&   4   &	0	&	922.03	&	104.96	&	1407.59	&	796.82	&	25.2	\\
            &   5   &	0	&	747.21	&	57.38	&	1073.74	&	687.80	&	71.0	\\
            &   5   &	6	&	742.94	&	15.60	&	908.80	&	706.06	&	95.2	\\
\vspace{0.5em}            &   5   &	12	&	761.68	&	10.00	&	789.06	&	726.57	&	100	\\
RC204C40	&   4   &	0	&	602.43	&	18.87	&	752.61	&	574.57	&	90.1	\\
            &   4   &	6	&	605.92	&	8.31	&	653.53	&	577.83	&	99.5	\\
            &   4   &	12	&	615.79	&	7.88	&	660.27	&	593.01	&	99.5	\\
\bottomrule
\end{tabular}
\label{table_simul_normal_correlated}
\end{table}

\begin{table}
\small
\renewcommand{\arraystretch}{0.8}
\centering 
\caption{Simulation Results - Normal Distribution}
\vspace{3mm}
\begin{tabular}{c c c c c c c c}
\toprule
\vspace{0.3em}
{$Inst.$} & {$n_v$} & {$\Gamma$} & {Avg.} & {Std.} & {Max} & {Min} & {Feasibility(\%)} \\
\hline
\rule{0pt}{2.7ex}
C103C40	&	5	&	0	&	807.6	&	138.81	&	1583.96	&	648.69	&	23.4	\\
	&	5	&	6	&	762.5	&	76.78	&	1155.82	&	671.94	&	47.5	\\
\vspace{0.5em}	&	5	&	12	&	723.73	&	36.21	&	1018.05	&	678.69	&	85.8	\\
C202C40	&	4	&	0	&	543.15	&	58.65	&	870.45	&	491.49	&	67.6	\\
    	&	4	&	6	&	538.24	&	20.94	&	710.38	&	505.69	&	91.4	\\
\vspace{0.5em}    	&	4	&	12	&	554.42	&	13.36	&	697.6	&	526.54	&	98.1	\\
C208C40	&	4	&	0	&	524.72	&	8.79	&	558.58	&	494.28	&	99.7	\\
    	&	4	&	6	&	525.38	&	8.85	&	552.52	&	492.43	&	99.9	\\
\vspace{0.5em}    	&	4	&	12	&	531.57	&	8.90	&	558.28	&	498.27	&	100	\\
R102C40	&	11	&	0	&	1181.63	&	207.84	&	2177.63	&	883.31	&	5.3	\\
	&	12	&	0	&	1138.62	&	178.29	&	2034.79	&	886.22	&	9.8	\\
	&	12	&	6	&	1036.07	&	112.16	&	1665.16	&	906.9	&	31.3	\\
\vspace{0.5em}	&	12	&	12	&	1013.62	&	96.04	&	1473.5	&	904.82	&	34.3	\\
R105C40	&	10	&	0	&	989.18	&	148.65	&	1705.83	&	779.64	&	10.2	\\
	&	11	&	0	&	995.98	&	152.66	&	1719.47	&	777.5	&	11.4	\\
	&	11	&	6	&	952.12	&	92.38	&	1466.38	&	845.67	&	39.7	\\
\vspace{0.5em}	&	11	&	12	&	953.44	&	84.52	&	1456.37	&	851.42	&	40.2	\\
R202C40	&	4	&	0	&	632.53	&	87.68	&	1068.21	&	543.23	&	43.1	\\
	&	4	&	6	&	606.42	&	52.02	&	942.59	&	556.72	&	70.8	\\
\vspace{0.5em}	&	4	&	12	&	597.37	&	34.08	&	835.21	&	550.43	&	81.6	\\
R209C40	&	4	&	0	&	508.47	&	33.24	&	766.05	&	475.85	&	88.5	\\
    	&	4	&	6	&	509.58	&	29.94	&	744.6	&	474.47	&	92.4	\\
\vspace{0.5em}    	&	4	&	12	&	527.35	&	19.69	&	766.51	&	494.48	&	94	\\
RC108C40	&	10	&	0	&	1067.62	&	172.8	&	2013.18	&	838.98	&	12.6	\\
	&	11	&	0	&	1129.68	&	187.01	&	2153.56	&	854.94	&	6.5	\\
	&	11	&	6	&	1018.84	&	117.12	&	1588.53	&	881.13	&	32.7	\\
\vspace{0.5em}	&	11	&	12	&	1077.39	&	122.02	&	1783.43	&	937.13	&	35.4	\\
RC202C40	&	4	&	0	&	1027.29	&	182.75	&	2045.43	&	786.35	&	12.3	\\
	&	5	&	0	&	792.86	&	115.78	&	1383.88	&	674.52	&	53.7	\\
	&	5	&	6	&	758.73	&	53.58	&	1215.99	&	694.5	&	79.3	\\
\vspace{0.5em}	&	5	&	12	&	766.68	&	24.19	&	963.26	&	719.13	&	91.8	\\
RC204C40	&	4	&	0	&	617.84	&	49.19	&	936.13	&	564.66	&	71.8	\\
	&	4	&	6	&	607.79	&	16.35	&	784.17	&	575.47	&	95.8	\\
	&	4	&	12	&	616.74	&	16.25	&	838.54	&	581.63	&	97.1	\\
\bottomrule
\end{tabular}
\label{table_simul_normal}
\end{table}

\begin{table}
\small
\renewcommand{\arraystretch}{0.8}
\centering 
\caption{Simulation Results - Uniform Distribution}
\vspace{3mm}
\begin{tabular}{c c c c c c c c}
\toprule
\vspace{0.3em}
{$Inst.$} & {$n_v$} & {$\Gamma$} & {Avg.} & {Std.} & {Max} & {Min} & {Feasibility(\%)} \\
\hline
\rule{0pt}{2.7ex}
C103C40	&   5   &	0	&	739.03	&	67.78	&	1036.19	&	665.46	&	35.4	\\
	    &   5   &	6	&	735.61	&	39.17	&	916.03	&	689.3	&	56.4	\\
\vspace{0.5em}	    &   5   &	12	&	714.83	&	7.08	&	768.01	&	692.08	&	99.8	\\
C202C40	&   4   &	0	&	523.72	&	22.90	&	676.86	&	501.50	&	80.9	\\
    	&   4   &	6	&	534.14	&	5.28	&	551.46	&	516.01	&	99.8	\\
\vspace{0.5em}    	&   4   &	12	&	553.49	&	5.32	&	572.52	&	537.16	&	100	\\
C208C40	&   4   &	0	&	524.81	&	5.06	&	539.40	&	511.50	&	100	\\
    	&   4   &	6	&	525.50	&	4.94	&	541.52	&	509.99	&	100	\\
\vspace{0.5em}    	&   4   &	12	&	531.12	&	5.06	&	547.28	&	516.28	&	100	\\
R102C40	&   11   &	0	&	1019.24	&	90.88	&	1390.84	&	893.74	&	12.2	\\
	    &   12    &	0	&	1018.97	&	84.91	&	1451.93	&	910.17	&	16.6	\\
	    &   12   &	6	&	981.84	&	52.09	&	1210.18	&	921.16	&	47.5	\\
\vspace{0.5em}	    &   12   &	12	&	959.03	&	34.7	&	1141.24	&	916.99	&	60.7	\\
R105C40	&   10   &	0	&	884.93	&	69.17	&	1213.68	&	799.75	&	21.8	\\
	    &   11   &	0	&	895.02	&	73.81	&	1231.86	&	804.69	&	22.7	\\
	    &   11   &	6	&	900.52	&	32.26	&	1083.69	&	866.46	&	72.2	\\
\vspace{0.5em}	    &   11   &	12	&	904.2	&	24.07	&	1047.31	&	873.99	&	77.4	\\
R202C40	&   4   &	0	&	598.56	&	43.49	&	799.6	&	552.54	&	51.9	\\
	    &   4   &	6	&	586.39	&	8.18	&	692.99	&	568.61	&	97.3	\\
\vspace{0.5em}	    &   4   &	12	&	587.54	&	7.65	&	660.87	&	569.93	&	97.7	\\
R209C40	&   4   &	0	&	500.09	&	6.42	&	587.13	&	484.23	&	99.5	\\
    	&   4   &	6	&	503.14	&	6.34	&	592.08	&	489.24	&	99.5	\\
\vspace{0.5em}    	&   4   &	12	&	523.94	&	5.18	&	542.28	&	505.40	&	99.9	\\
RC108C40	&   10   &	0	&	963.13	&	88.27	&	1337.97	&	847.89	&	18.7	\\
	    &   11   &	0	&	1009.42	&	97.87	&	1524.18	&	864.14	&	12	\\
	    &   11   &	6	&	946.31	&	31.76	&	1163.8	&	907.45	&	70.9	\\
\vspace{0.5em}	    &   11   &	12	&	998.57	&	32.13	&	1153.84	&	954.78	&	70.3	\\
RC202C40	&   4   &	0	&	912.81	&	93.49	&	1320.14	&	805.52	&	25.3	\\
	        &   5   &	0	&	746.75	&	51.55	&	1034.57	&	700.84	&	70	\\
	        &   5   &	6	&	742.9	&	11.6	&	836.07	&	716.96	&	95.1	\\
\vspace{0.5em}	        &   5   &	12	&	761.85	&	7.09	&	782.61	&	737.85	&	100	\\
RC204C40	&   4   &	0	&	600.61	&	12.34	&	714.16	&	580.96	&	92.3	\\
	        &   4   &	6	&	605.94	&	5.87	&	625.34	&	589.74	&	99.9	\\
	        &   4   &	12	&	615.61	&	5.67	&	631.65	&	598.37	&	100	\\
\bottomrule
\end{tabular}
\label{table_simul_uniform}
\end{table}

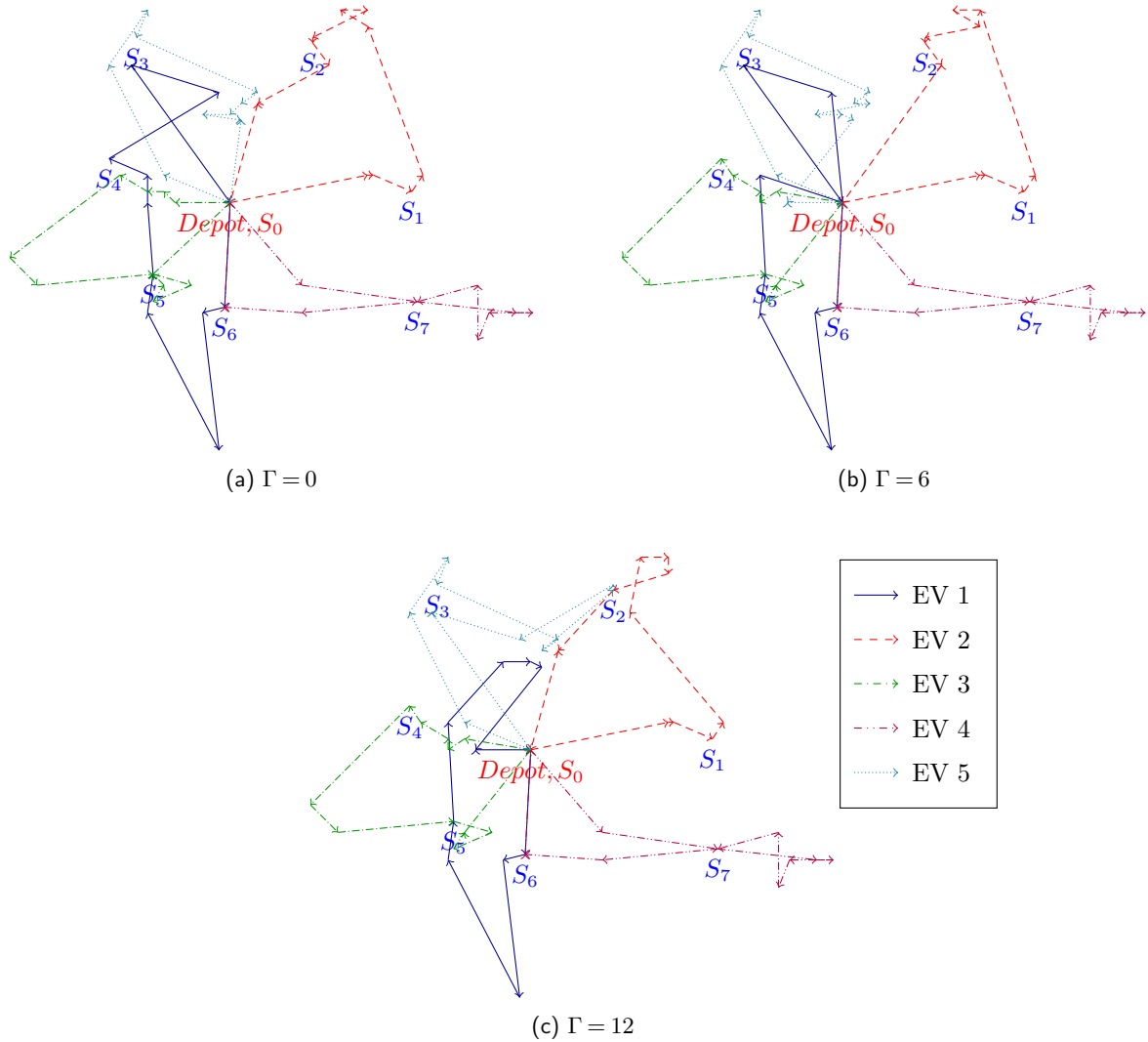
\begin{figure}[!ht]
     \centering
     \begin{subfigure}[b]{0.45\textwidth}
         \centering
         \begin{tikzpicture}[scale=0.075]
            \def\D{\textcolor{red}{\small $Depot, S_0$}}
            \def\SA{\textcolor{blue}{\small $S_1$}}
            \def\SB{\textcolor{blue}{\small $S_2$}}
            \def\SC{\textcolor{blue}{\small $S_3$}}
            \def\SD{\textcolor{blue}{\small $S_4$}}
            \def\SE{\textcolor{blue}{\small $S_5$}}
            \def\SF{\textcolor{blue}{\small $S_6$}}
            \def\SG{\textcolor{blue}{\small $S_7$}}
            \colorlet{input}{blue!60!black}
            \colorlet{route1}{red}
            \colorlet{route2}{green!60!black}
            \colorlet{route3}{purple}
            \colorlet{route4}{cyan!60!black}
            \coordinate [label=below:\D] (0) at (40, 50);
            \coordinate (1) at (22, 75);
            \coordinate (2) at (45, 68);
            \coordinate (3) at (58, 75);
            \coordinate (4) at (42, 65);
            \coordinate (5) at (18, 75);
            \coordinate (6) at (25, 30);
            \coordinate (7) at (33, 35);
            \coordinate (8) at (35, 66);
            \coordinate (9) at (35, 30);
            \coordinate (10) at (85, 25);
            \coordinate (11) at (28, 35);
            \coordinate (12) at (65, 55);
            \coordinate (13) at (38, 70);
            \coordinate (14) at (53, 30);
            \coordinate (15) at (28, 55);
            \coordinate (16) at (66, 55);
            \coordinate (17) at (53, 35);
            \coordinate (18) at (25, 55);
            \coordinate (19) at (65, 82);
            \coordinate (20) at (30, 50);
            \coordinate (21) at (0, 40);
            \coordinate (22) at (95, 30);
            \coordinate (23) at (28, 52);
            \coordinate (24) at (26, 32);
            \coordinate (25) at (25, 50);
            \coordinate (26) at (40, 66);
            \coordinate (27) at (45, 70);
            \coordinate (28) at (60, 85);
            \coordinate (29) at (75, 55);
            \coordinate (30) at (20, 55);
            \coordinate (31) at (65, 85);
            \coordinate (32) at (92, 30);
            \coordinate (33) at (25, 52);
            \coordinate (34) at (25, 85);
            \coordinate (35) at (5, 35);
            \coordinate (36) at (85, 35);
            \coordinate (37) at (38, 5);
            \coordinate (38) at (55, 80);
            \coordinate (39) at (42, 68);
            \coordinate (40) at (87, 30);
            \coordinate (41) at (40, 50);
            \coordinate[label=below:\SA] (42) at (73, 52);
            \coordinate[label=below:\SB] (43) at (55, 79);
            \coordinate[label=below:\SC] (44) at (23, 80);
            \coordinate[label=below:\SD] (45) at (18, 58);
            \coordinate[label=below:\SE] (46) at (26, 37);
            \coordinate[label=below:\SF] (47) at (39, 31);
            \coordinate[label=below:\SG] (48) at (74, 32);
            \coordinate (49) at (40, 50);
            \coordinate (50) at (73, 52);
            \coordinate (51) at (55, 79);
            \coordinate (52) at (23, 80);
            \coordinate (53) at (18, 58);
            \coordinate (54) at (26, 37);
            \coordinate (55) at (39, 31);
            \coordinate (56) at (74, 32);
            \coordinate (57) at (40, 50);
            \draw [->, color=input, solid] (0) edge (47) (47) edge (9) (9) edge (37) (37) edge (6) (6) edge (46) (46) edge (25) (25) edge (18) (18) edge (45) (45) edge (13) (13) edge (1) (1) edge (57);
            \draw [->, color=route1, densely dashed] (0) edge (12) (12) edge (16) (16) edge (42) (42) edge (29) (29) edge (19) (19) edge (28) (28) edge (31) (31) edge (38) (38) edge (43) (43) edge (3) (3) edge (2) (2) edge (57);
            \draw [->, color=route2, densely dashdotted] (0) edge (20) (20) edge (23) (23) edge (33) (33) edge (30) (30) edge (21) (21) edge (35) (35) edge (46) (46) edge (7) (7) edge (24) (24) edge (11) (11) edge (54) (54) edge (57);
            \draw [->, color=route3, densely dashdotdotted] (0) edge (17) (17) edge (56) (56) edge (32) (32) edge (22) (22) edge (40) (40) edge (10) (10) edge (36) (36) edge (48) (48) edge (14) (14) edge (47) (47) edge (57);
            \draw [->, color=route4, densely dotted] (0) edge (15) (15) edge (5) (5) edge (34) (34) edge (44) (44) edge (27) (27) edge (39) (39) edge (26) (26) edge (8) (8) edge (4) (4) edge (57);
            % Deterministic : [[0, 47, 9, 37, 6, 46, 25, 18, 45, 13, 1, 57], [0, 12, 16, 42, 29, 19, 28, 31, 38, 43, 3, 2, 57], [0, 20, 23, 33, 30, 21, 35, 46, 7, 24, 11, 54, 57], [0, 17, 56, 32, 22, 40, 10, 36, 48, 14, 47, 57], [0, 15, 5, 34, 44, 27, 39, 26, 8, 4, 57]]
        
        \end{tikzpicture}
         \caption{$\Gamma=0$}
         \label{fig_budget_0_c103}
     \end{subfigure}
     \hspace{1.5em}
     \begin{subfigure}[b]{0.45\textwidth}
         \centering
         \begin{tikzpicture}[scale=0.075]
            \def\D{\textcolor{red}{\small $Depot, S_0$}}
            \def\SA{\textcolor{blue}{\small $S_1$}}
            \def\SB{\textcolor{blue}{\small $S_2$}}
            \def\SC{\textcolor{blue}{\small $S_3$}}
            \def\SD{\textcolor{blue}{\small $S_4$}}
            \def\SE{\textcolor{blue}{\small $S_5$}}
            \def\SF{\textcolor{blue}{\small $S_6$}}
            \def\SG{\textcolor{blue}{\small $S_7$}}
            \colorlet{input}{blue!60!black}
            \colorlet{route1}{red}
            \colorlet{route2}{green!60!black}
            \colorlet{route3}{purple}
            \colorlet{route4}{cyan!60!black}
            \coordinate [label=below:\D] (0) at (40, 50);
            \coordinate (1) at (22, 75);
            \coordinate (2) at (45, 68);
            \coordinate (3) at (58, 75);
            \coordinate (4) at (42, 65);
            \coordinate (5) at (18, 75);
            \coordinate (6) at (25, 30);
            \coordinate (7) at (33, 35);
            \coordinate (8) at (35, 66);
            \coordinate (9) at (35, 30);
            \coordinate (10) at (85, 25);
            \coordinate (11) at (28, 35);
            \coordinate (12) at (65, 55);
            \coordinate (13) at (38, 70);
            \coordinate (14) at (53, 30);
            \coordinate (15) at (28, 55);
            \coordinate (16) at (66, 55);
            \coordinate (17) at (53, 35);
            \coordinate (18) at (25, 55);
            \coordinate (19) at (65, 82);
            \coordinate (20) at (30, 50);
            \coordinate (21) at (0, 40);
            \coordinate (22) at (95, 30);
            \coordinate (23) at (28, 52);
            \coordinate (24) at (26, 32);
            \coordinate (25) at (25, 50);
            \coordinate (26) at (40, 66);
            \coordinate (27) at (45, 70);
            \coordinate (28) at (60, 85);
            \coordinate (29) at (75, 55);
            \coordinate (30) at (20, 55);
            \coordinate (31) at (65, 85);
            \coordinate (32) at (92, 30);
            \coordinate (33) at (25, 52);
            \coordinate (34) at (25, 85);
            \coordinate (35) at (5, 35);
            \coordinate (36) at (85, 35);
            \coordinate (37) at (38, 5);
            \coordinate (38) at (55, 80);
            \coordinate (39) at (42, 68);
            \coordinate (40) at (87, 30);
            \coordinate (41) at (40, 50);
            \coordinate[label=below:\SA] (42) at (73, 52);
            \coordinate[label=below:\SB] (43) at (55, 79);
            \coordinate[label=below:\SC] (44) at (23, 80);
            \coordinate[label=below:\SD] (45) at (18, 58);
            \coordinate[label=below:\SE] (46) at (26, 37);
            \coordinate[label=below:\SF] (47) at (39, 31);
            \coordinate[label=below:\SG] (48) at (74, 32);
            \coordinate (49) at (40, 50);
            \coordinate (50) at (73, 52);
            \coordinate (51) at (55, 79);
            \coordinate (52) at (23, 80);
            \coordinate (53) at (18, 58);
            \coordinate (54) at (26, 37);
            \coordinate (55) at (39, 31);
            \coordinate (56) at (74, 32);
            \coordinate (57) at (40, 50);
            \draw [->, color=input, solid] (0) edge (47) (47) edge (9) (9) edge (37) (37) edge (6) (6) edge (54) (54) edge (18) (18) edge (41) (41) edge (13) (13) edge (1) (1) edge (57);
            \draw [->, color=route1, densely dashed] (0) edge (12) (12) edge (16) (16) edge (42) (42) edge (29) (29) edge (31) (31) edge (28) (28) edge (19) (19) edge (38) (38) edge (43) (43) edge (3) (3) edge (57);
            \draw [->, color=route2, densely dashdotted] (0) edge (23) (23) edge (25) (25) edge (33) (33) edge (30) (30) edge (45) (45) edge (21) (21) edge (35) (35) edge (46) (46) edge (7) (7) edge (24) (24) edge (11) (11) edge (57);
            \draw [->, color=route3, densely dashdotdotted] (0) edge (17) (17) edge (48) (48) edge (32) (32) edge (22) (22) edge (40) (40) edge (10) (10) edge (36) (36) edge (56) (56) edge (14) (14) edge (47) (47) edge (57);
            \draw [->, color=route4, densely dotted] (0) edge (15) (15) edge (5) (5) edge (34) (34) edge (44) (44) edge (27) (27) edge (39) (39) edge (2) (2) edge (8) (8) edge (26) (26) edge (4) (4) edge (20) (20) edge (57);
            % Budget = 6 : [[0, 47, 9, 37, 6, 54, 18, 41, 13, 1, 57], [0, 15, 5, 34, 44, 27, 39, 2, 8, 26, 4, 20, 57], [0, 17, 48, 32, 22, 40, 10, 36, 56, 14, 47, 57], [0, 23, 25, 33, 30, 45, 21, 35, 46, 7, 24, 11, 57], [0, 12, 16, 42, 29, 31, 28, 19, 38, 43, 3, 57]]
        
        \end{tikzpicture}
         \caption{$\Gamma = 6$}
         \label{fig_budget_6_c103}
     \end{subfigure}
     \vspace{2em}
     \\
     \begin{subfigure}[b]{0.45\textwidth}
        \centering
        \begin{tikzpicture}[scale=0.075,
        R1/.style={shape=circle, draw=black, line width=2},
        R2/.style={shape=circle, draw=blue, line width=2},
        R3/.style={shape=circle, draw=green, line width=2},
        R4/.style={shape=circle, draw=red, line width=2}]
            \def\D{\textcolor{red}{\small $Depot, S_0$}}
            \def\SA{\textcolor{blue}{\small $S_1$}}
            \def\SB{\textcolor{blue}{\small $S_2$}}
            \def\SC{\textcolor{blue}{\small $S_3$}}
            \def\SD{\textcolor{blue}{\small $S_4$}}
            \def\SE{\textcolor{blue}{\small $S_5$}}
            \def\SF{\textcolor{blue}{\small $S_6$}}
            \def\SG{\textcolor{blue}{\small $S_7$}}
            \colorlet{input}{blue!60!black}
            \colorlet{route1}{red}
            \colorlet{route2}{green!60!black}
            \colorlet{route3}{purple}
            \colorlet{route4}{cyan!60!black}
            \coordinate [label=below:\D] (0) at (40, 50);
            \coordinate (1) at (22, 75);
            \coordinate (2) at (45, 68);
            \coordinate (3) at (58, 75);
            \coordinate (4) at (42, 65);
            \coordinate (5) at (18, 75);
            \coordinate (6) at (25, 30);
            \coordinate (7) at (33, 35);
            \coordinate (8) at (35, 66);
            \coordinate (9) at (35, 30);
            \coordinate (10) at (85, 25);
            \coordinate (11) at (28, 35);
            \coordinate (12) at (65, 55);
            \coordinate (13) at (38, 70);
            \coordinate (14) at (53, 30);
            \coordinate (15) at (28, 55);
            \coordinate (16) at (66, 55);
            \coordinate (17) at (53, 35);
            \coordinate (18) at (25, 55);
            \coordinate (19) at (65, 82);
            \coordinate (20) at (30, 50);
            \coordinate (21) at (0, 40);
            \coordinate (22) at (95, 30);
            \coordinate (23) at (28, 52);
            \coordinate (24) at (26, 32);
            \coordinate (25) at (25, 50);
            \coordinate (26) at (40, 66);
            \coordinate (27) at (45, 70);
            \coordinate (28) at (60, 85);
            \coordinate (29) at (75, 55);
            \coordinate (30) at (20, 55);
            \coordinate (31) at (65, 85);
            \coordinate (32) at (92, 30);
            \coordinate (33) at (25, 52);
            \coordinate (34) at (25, 85);
            \coordinate (35) at (5, 35);
            \coordinate (36) at (85, 35);
            \coordinate (37) at (38, 5);
            \coordinate (38) at (55, 80);
            \coordinate (39) at (42, 68);
            \coordinate (40) at (87, 30);
            \coordinate (41) at (40, 50);
            \coordinate[label=below:\SA] (42) at (73, 52);
            \coordinate[label=below:\SB] (43) at (55, 79);
            \coordinate[label=below:\SC] (44) at (23, 80);
            \coordinate[label=below:\SD] (45) at (18, 58);
            \coordinate[label=below:\SE] (46) at (26, 37);
            \coordinate[label=below:\SF] (47) at (39, 31);
            \coordinate[label=below:\SG] (48) at (74, 32);
            \coordinate (49) at (40, 50);
            \coordinate (50) at (73, 52);
            \coordinate (51) at (55, 79);
            \coordinate (52) at (23, 80);
            \coordinate (53) at (18, 58);
            \coordinate (54) at (26, 37);
            \coordinate (55) at (39, 31);
            \coordinate (56) at (74, 32);
            \coordinate (57) at (40, 50);
            \draw [->, color=input, solid] (0) edge (47) (47) edge (9) (9) edge (37) (37) edge (6) (6) edge (46) (46) edge (18) (18) edge (8) (8) edge (26) (26) edge (4) (4) edge (20) (20) edge (57);
            \draw [->, color=route1, densely dashed] (0) edge (12) (12) edge (16) (16) edge (42) (42) edge (29) (29) edge (3) (3) edge (28) (28) edge (31) (31) edge (19) (19) edge (51) (51) edge (2) (2) edge (57);
            \draw [->, color=route2, densely dashdotted] (0) edge (23) (23) edge (25) (25) edge (33) (33) edge (30) (30) edge (45) (45) edge (21) (21) edge (35) (35) edge (54) (54) edge (7) (7) edge (24) (24) edge (11) (11) edge (57);
            \draw [->, color=route3, densely dashdotdotted] (0) edge (17) (17) edge (48) (48) edge (32) (32) edge (22) (22) edge (40) (40) edge (10) (10) edge (36) (36) edge (56) (56) edge (14) (14) edge (47) (47) edge (57);
            \draw [->, color=route4, densely dotted] (0) edge (15) (15) edge (5) (5) edge (34) (34) edge (44) (44) edge (27) (27) edge (39) (39) edge (43) (43) edge (38) (38) edge (13) (13) edge (1) (1) edge (57);
            
            % Budget = 12 : [[0, 15, 5, 34, 44, 27, 39, 43, 38, 13, 1, 57], [0, 12, 16, 42, 29, 3, 28, 31, 19, 51, 2, 57], [0, 17, 48, 32, 22, 40, 10, 36, 56, 14, 47, 57], [0, 23, 25, 33, 30, 45, 21, 35, 54, 7, 24, 11, 57], [0, 47, 9, 37, 6, 46, 18, 8, 26, 4, 20, 57]]
        
        \path ([xshift=1cm,yshift=-5mm]current bounding box.north east)
         node[matrix,anchor=north west,cells={nodes={anchor=west}},
         draw,inner sep=1ex]{
          \draw[-to,color=input,solid](0,0) -- ++ (0.6,0); & \node{\small EV 1};\\ % Solid line
          \draw[-to,color=route1,dashed](0,0) -- ++ (0.6,0); & \node{\small EV 2};\\ %Dashed line
          \draw[-to,color=route2,dashdotted](0,0) -- ++ (0.6,0); & \node{\small EV 3};\\ %Alternating dashed and dotted line
          \draw[-to,color=route3,dashdotdotted](0,0) -- ++ (0.6,0); & \node{\small EV 4};\\ %Alternating dashed and double-dotted line
          \draw[-to,color=route4,densely dotted](0,0) -- ++ (0.6,0); & \node{\small EV 5};\\ %Densely dotted line
         };
        
        \end{tikzpicture}
         \caption{$\Gamma = 12$}
         \label{fig_budget_12_c103}
     \end{subfigure}
        \caption{EV Routes - Instance of C103C40}
        \label{fig_routes_visualization}
\end{figure}

To provide more managerial insights, we choose one mid-sized instance (C103C40) and illustrate the route of each model in terms of the amount of electricity charging en route, travel ending times, and idle times. We first visualize the routes of the five EVs for three cases in Figure \ref{fig_routes_visualization}. Figure \ref{fig_budget_0_c103} represents the EVs' routes of the deterministic model ($\Gamma = 0$), and Figures \ref{fig_budget_6_c103} and \ref{fig_budget_12_c103} represent the routes of ARO models with $\Gamma =  6$ and $\Gamma = 12$, respectively. We can observe that each EV visits slightly different customers and stations for each model. In all three models, the total number of visited stations is the same. However, the number of visited customers for each EV can be different for each model. For instance, for the deterministic model, the five EVs visit 7, 9, 9, 7, and 8 customers in order of EV1 to EV5. On the other hand, for the ARO models with $\Gamma =  6$ the number of customers for each EV is: 6, 8, 9, 7, and 10 customers and for $\Gamma =  12$ there are: 8, 8, 9, 7, and 8 customers.

\begin{figure}[!ht]
     \centering
     \begin{subfigure}[b]{0.35\textwidth}
         \centering
        \begin{tikzpicture}
        \begin{axis}[xtick={80,90,100,110,120}, xticklabels={80,90,100,110,120}, xmin=73, xmax=125, ytick={1,2,3}, yticklabels={0,6,12}, ylabel={$\Gamma$}, ylabel style={rotate=-90}, width=2.4in, height=2.1in, /pgfplots/boxplot/box extend=0.4,
        label style={font=\small},
        tick label style={font=\small}]
        \addplot+ [boxplot prepared={lower whisker=75.536, lower quartile=87.949, median=90.724, upper quartile=96.362, upper whisker=105.451},] table[row sep=\\,y index=0] {\\};
        \addplot+ [boxplot prepared={lower whisker=83.931, lower quartile=91.29, median=94.37, upper quartile=97.08, upper whisker=102.57},] table[row sep=\\,y index=0] {\\};
        \addplot+ [mark=*, boxplot prepared={lower whisker=88.56, lower quartile=100.51, median=105.73, upper quartile=108.56, upper whisker=114.31},] table[row sep=\\,y index=0] {88.4\\};
        \end{axis}
        \end{tikzpicture}
         \caption{\footnotesize Normal distribution with correlation}
         \label{boxplot_charge_normal_correlated}
     \end{subfigure}
     \begin{subfigure}[b]{0.32\textwidth}
         \centering
         \begin{tikzpicture}
        \begin{axis}[xtick={80,90,100,110,120}, xticklabels={80,90,100,110,120}, xmin=73, xmax=125, ytick={1,2,3}, yticklabels={0,6,12}, %ylabel={$\Gamma$}, 
        ylabel style={rotate=-90}, width=2.4in, height=2.1in, /pgfplots/boxplot/box extend=0.4,
        label style={font=\small},
        tick label style={font=\small}]
        \addplot+ [mark=*, boxplot prepared={lower whisker=77.749, lower quartile=84.14, median=89.39, upper quartile=93.52, upper whisker=104.32},] table[row sep=\\,y index=0] {\\};
        \addplot+ [mark=*, boxplot prepared={lower whisker=83.2, lower quartile=91.29, median=94.37, upper quartile=96.94, upper whisker=105.19},] table[row sep=\\,y index=0] {82.69\\};
        \addplot+ [boxplot prepared={lower whisker=83.55, lower quartile=96.92, median=103.88, upper quartile=108.35, upper whisker=115.88},] table[row sep=\\,y index=0] {\\};
        \end{axis}
        \end{tikzpicture}
         \caption{\footnotesize Normal distribution}
         \label{boxplot_charge_normal}
     \end{subfigure}
     \begin{subfigure}[b]{0.31\textwidth}
        \centering
        \begin{tikzpicture}
        \begin{axis}[xtick={80,90,100,110,120}, xticklabels={80,90,100,110,120}, xmin=73, xmax=125, ytick={1,2,3}, yticklabels={0,6,12}, %ylabel={$\Gamma$}, 
        ylabel style={rotate=-90}, width=2.4in, height=2.1in, /pgfplots/boxplot/box extend=0.4,
        label style={font=\small},
        tick label style={font=\small}]
        \addplot+ [mark=*, boxplot prepared={lower whisker=80.14, lower quartile=89.69, median=92.2, upper quartile=100.7, upper whisker=105.97},] table[row sep=\\,y index=0] {\\};
        \addplot+ [mark=*, boxplot prepared={lower whisker=86.62, lower quartile=92.68, median=94.75, upper quartile=96.88, upper whisker=102.84},] table[row sep=\\,y index=0] {103.54\\};
        \addplot+ [mark=*, boxplot prepared={lower whisker=88.59, lower quartile=98.96, median=106.64, upper quartile=108.83, upper whisker=115.67},] table[row sep=\\,y index=0] {\\};
        \end{axis}
        \end{tikzpicture}
         \caption{\footnotesize Uniform distribution}
         \label{boxplot_charge_uniform}
     \end{subfigure}
        \caption{Total recharging amount(MWh)}
        \label{fig_boxplot_charge}
\end{figure}
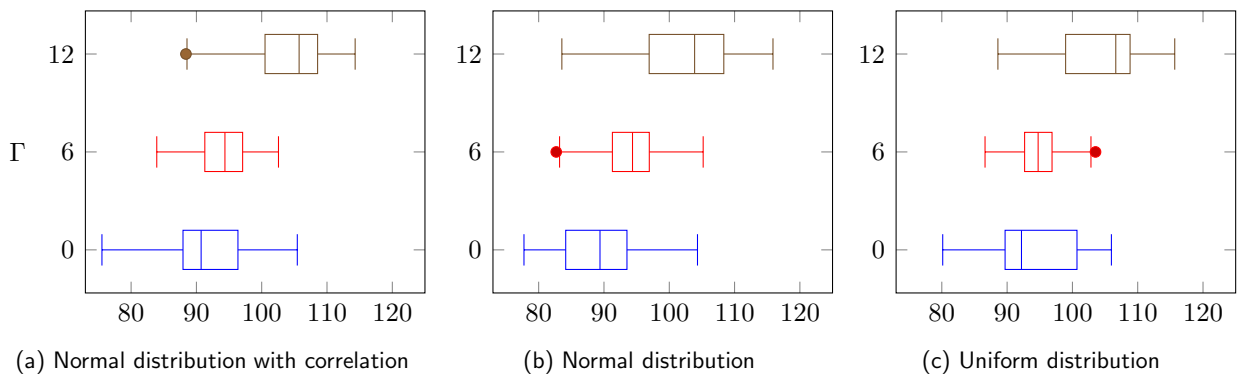

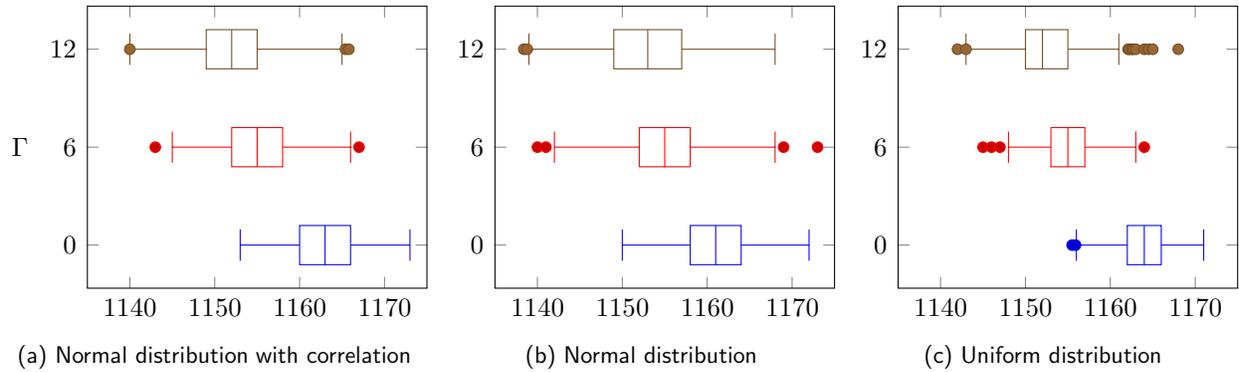
\begin{figure}[!ht]
     \centering
     \begin{subfigure}[b]{0.35\textwidth}
         \centering
        \begin{tikzpicture}
        \begin{axis}[xtick={1140,1150,1160,1170}, xticklabels={1140,1150,1160,1170}, xmin=1135, xmax=1175, ytick={1,2,3}, yticklabels={0,6,12}, ylabel={$\Gamma$}, ylabel style={rotate=-90}, width=2.4in, height=2.1in, /pgfplots/boxplot/box extend=0.4,
        label style={font=\small},
        tick label style={font=\small}]
        \addplot+ [mark=*, boxplot prepared={lower whisker=1153, lower quartile=1160, median=1163, upper quartile=1166, upper whisker=1173},] table[row sep=\\,y index=0] {\\};
        \addplot+ [mark=*, boxplot prepared={lower whisker=1145, lower quartile=1152, median=1155, upper quartile=1158, upper whisker=1166},] table[row sep=\\,y index=0] {1143\\1167\\};
        \addplot+ [mark=*, boxplot prepared={lower whisker=1140, lower quartile=1149, median=1152, upper quartile=1155, upper whisker=1165},] table[row sep=\\,y index=0] {1140\\1165.4\\1165.8\\};
        \end{axis}
        \end{tikzpicture}
         \caption{\footnotesize Normal distribution with correlation}
         \label{boxplot_end_normal_correlated}
     \end{subfigure}
     \begin{subfigure}[b]{0.32\textwidth}
         \centering
         \begin{tikzpicture}
        \begin{axis}[xtick={1140,1150,1160,1170}, xticklabels={1140,1150,1160,1170}, xmin=1135, xmax=1175, ytick={1,2,3}, yticklabels={0,6,12}, %ylabel={$\Gamma$}, 
        ylabel style={rotate=-90}, width=2.4in, height=2.1in, /pgfplots/boxplot/box extend=0.4,
        label style={font=\small},
        tick label style={font=\small}]
        \addplot+ [mark=*, boxplot prepared={lower whisker=1150, lower quartile=1158, median=1161, upper quartile=1164, upper whisker=1172},] table[row sep=\\,y index=0] {\\};
        \addplot+ [mark=*, boxplot prepared={lower whisker=1142, lower quartile=1152, median=1155, upper quartile=1158, upper whisker=1168},] table[row sep=\\,y index=0] {1140\\1141\\1169\\1173\\};
        \addplot+ [mark=*, boxplot prepared={lower whisker=1139, lower quartile=1149, median=1153, upper quartile=1157, upper whisker=1168},] table[row sep=\\,y index=0] {1138.4\\1138.8\\};
        \end{axis}
        \end{tikzpicture}
         \caption{\footnotesize Normal distribution}
         \label{boxplot_end_normal}
     \end{subfigure}
     \begin{subfigure}[b]{0.31\textwidth}
        \centering
        \begin{tikzpicture}
        \begin{axis}[xtick={1140,1150,1160,1170}, xticklabels={1140,1150,1160,1170}, xmin=1135, xmax=1175, ytick={1,2,3}, yticklabels={0,6,12}, %ylabel={$\Gamma$}, 
        ylabel style={rotate=-90}, width=2.4in, height=2.1in, /pgfplots/boxplot/box extend=0.4,
        label style={font=\small},
        tick label style={font=\small}]
        \addplot+ [mark=*, boxplot prepared={lower whisker=1156, lower quartile=1162, median=1164, upper quartile=1166, upper whisker=1171},] table[row sep=\\,y index=0] {1155.9\\1155.5\\};
        \addplot+ [mark=*, boxplot prepared={lower whisker=1148, lower quartile=1153, median=1155, upper quartile=1157, upper whisker=1163},] table[row sep=\\,y index=0] {1145\\1146\\1147\\1164\\};
        \addplot+ [mark=*, boxplot prepared={lower whisker=1143, lower quartile=1150, median=1152, upper quartile=1155, upper whisker=1161},] table[row sep=\\,y index=0] {1142\\1143\\1162.1\\1162.2\\1162.5\\1162.7\\1163\\1164\\1164.5\\1165\\1168\\};
        \end{axis}
        \end{tikzpicture}
         \caption{\footnotesize Uniform distribution}
         \label{boxplot_end_uniform}
     \end{subfigure}
        \caption{Average travel ending time (min.)}
        \label{fig_boxplot_end}
\end{figure}

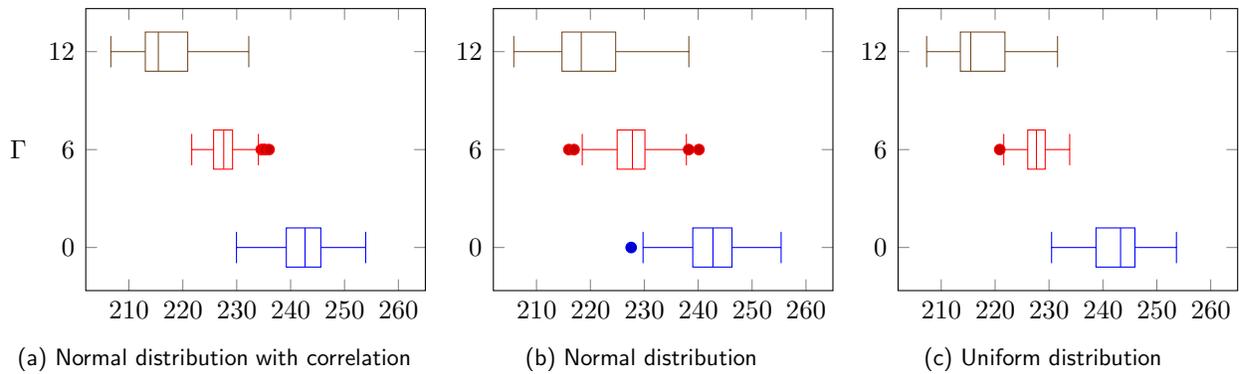
\begin{figure}[!ht]
     \centering
     \begin{subfigure}[b]{0.35\textwidth}
        \centering
        \begin{tikzpicture}
        \begin{axis}[xtick={210,220,230,240,250,260}, xticklabels={210,220,230,240,250,260}, xmin=202, xmax=265, ytick={1,2,3}, yticklabels={0,6,12}, ylabel={$\Gamma$}, ylabel style={rotate=-90}, width=2.4in, height=2.1in, /pgfplots/boxplot/box extend=0.4,
        label style={font=\small},
        tick label style={font=\small}]
        \addplot+ [boxplot prepared={lower whisker=229.958, lower quartile=239.203, median=242.68, upper quartile=245.604, upper whisker=253.935},] table[row sep=\\,y index=0] {\\};
        \addplot+ [mark=*, boxplot prepared={lower whisker=221.655, lower quartile=225.713, median=227.593, upper quartile=229.231, upper whisker=234.025},] table[row sep=\\,y index=0] {234.56\\234.97\\235.3\\236\\};
        \addplot+ [boxplot prepared={lower whisker=206.651, lower quartile=213.044, median=215.46, upper quartile=220.91, upper whisker=232.24},] table[row sep=\\,y index=0] {\\};
        \end{axis}
        \end{tikzpicture}
         \caption{\footnotesize Normal distribution with correlation}
         \label{boxplot_normal_correlated}
     \end{subfigure}
     \begin{subfigure}[b]{0.32\textwidth}
        \centering
        \begin{tikzpicture}
        \begin{axis}[xtick={210,220,230,240,250,260}, xticklabels={210,220,230,240,250,260}, xmin=202, xmax=265, ytick={1,2,3},yticklabels={0,6,12}, %ylabel={$\Gamma$},
        ylabel style={rotate=-90}, width=2.4in, height=2.1in, /pgfplots/boxplot/box extend=0.4,
        label style={font=\small},
        tick label style={font=\small}]
        \addplot+ [mark=*, boxplot prepared={lower whisker=229.7729, lower quartile=239, median=242.748, upper quartile=246.275, upper whisker=255.387},] table[row sep=\\,y index=0] {227.566\\};
        \addplot+ [mark=*, boxplot prepared={lower whisker=218.49, lower quartile=224.97, median=227.8, upper quartile=230.12, upper whisker=237.81},] table[row sep=\\,y index=0] {216.006\\216.97\\238.23\\240.14\\};
        \addplot+ [boxplot prepared={lower whisker=205.79, lower quartile=214.73, median=218.31, upper quartile=224.69, upper whisker=238.27},] table[row sep=\\,y index=0] {\\};
        \end{axis}
        \end{tikzpicture}
         \caption{\footnotesize Normal distribution}
         \label{boxplot_normal}
     \end{subfigure}
     \begin{subfigure}[b]{0.31\textwidth}
        \centering
        \begin{tikzpicture}
        \begin{axis}[xtick={210,220,230,240,250,260}, xticklabels={210,220,230,240,250,260}, xmin=202, xmax=265, ytick={1,2,3},yticklabels={0,6,12},%ylabel={$\Gamma$},
        ylabel style={rotate=-90}, width=2.4in, height=2.1in, /pgfplots/boxplot/box extend=0.4,
        label style={font=\small},
        tick label style={font=\small}]
        \addplot+ [mark=*, boxplot prepared={lower whisker=230.479, lower quartile=238.717, median=243.277, upper quartile=245.937, upper whisker=253.641},] table[row sep=\\,y index=0] {\\};
        \addplot+ [mark=*, boxplot prepared={lower whisker=221.6067, lower quartile=226.052, median=227.67, upper quartile=229.295, upper whisker=233.84},] table[row sep=\\,y index=0] {220.86\\};
        \addplot+ [mark=*, boxplot prepared={lower whisker=207.297, lower quartile=213.554, median=215.479, upper quartile=221.824, upper whisker=231.58},] table[row sep=\\,y index=0] {\\};
        \end{axis}
        \end{tikzpicture}
         \caption{\footnotesize Uniform distribution}
         \label{boxplot_uniform}
     \end{subfigure}
        \caption{Average idle time (min.)}
        \label{fig_boxplot_idle}
\end{figure}

Figures \ref{fig_boxplot_charge}, \ref{fig_boxplot_end}, and \ref{fig_boxplot_idle} illustrate box plots of the total recharging amount en routes, the average travel ending time, and the average idle time respectively. We use the data from the simulation tests presented in Tables \ref{table_simul_normal_correlated}, \ref{table_simul_normal}, and \ref{table_simul_uniform}. For each model, only the feasible scenarios are used for the box plot. Figure \ref{fig_boxplot_charge} shows that the total recharging amount en route increases as Gamma increases. Since we report only feasible scenarios, Figure \ref{fig_boxplot_charge} can also be interpreted as the total recharging amount. It means that routes from ARO models have more time flexibility for recharging than the deterministic ones. Thus, ARO models show higher feasibility than the deterministic model (see Table \ref{table_simul_normal_correlated}, \ref{table_simul_normal}, \ref{table_simul_uniform}). Figure \ref{fig_boxplot_end} shows that the average end-travel time of EVs decreases as $\Gamma$ increases. It denotes that routes of ARO models spend fewer time on the road than the deterministic model. Figure \ref{fig_boxplot_idle} shows the average idle time of EVs. Average idle time also tends to decrease as $\Gamma$ increases. With the results of Figures \ref{fig_boxplot_end}, and \ref{fig_boxplot_idle}, we can demonstrate that the ARO model can provide a route that utilizes EV fleets more efficiently. More precisely, a route from the ARO model can charge more electricity en route, may spend less time waiting for the next visit, and may complete the service earlier than the deterministic model. 

\section{Conclusion}\label{conclusion}

This paper presents an adaptive robust model for the electric vehicle routing problem with time window and partial recharging when the energy consumption rate exhibits uncertainty. To explain the necessity of an adaptive robust optimization approach, we provide a robust model and show that robust optimization approach might not be appropriate for the EVRPTWPR.  We postulate that recharging decision can be made after uncertainty realization to complete an EV's travel without any time window violation or negative battery level. Since recharging decisions can affect other decisions, such as service start time and battery level at each vertex, the proposed model is formulated as a two-stage adaptive robust problem. A solution method based on the column-and-constraint generation framework has been proposed to solve our model in a reasonable computational time.

Small and mid-sized instances are used to demonstrate the advantage of the proposed model. Numerical results show the economic efficiency and robustness of the proposed model. Robust routes from our model may require more electrical energy compared to the deterministic model. However, they are more secure and reliable than the routes given by the deterministic model. In the adaptive robust model, there is a tradeoff between the expected amount of electricity charged and the robustness of routes. Decision makers can decide the proper tradeoff they need to consider to optimize their operations. 

For future research, the {exact method} for the proposed model can be investigated. Since the proposed solution method is based on several heuristic algorithms, it faces some numerical instability. Also, extensive sensitivity analysis for the adaptive robust model with various uncertainty sets can be considered. As we mentioned in Section \ref{section_simulation}, uncertainties in the real world could be highly correlated. Furthermore, depending on the structure of the uncertainty set, there could be other proposed approaches to consider. Lastly, the time-variant uncertainty of energy consumption rate can be considered as future work. For example, when different EVs travel the same arcs but at different hours, each EV may have a different energy consumption rate. Therefore, introducing a time factor has practical significance since it could reflect more actual operations.

% Acknowledgments here
% \Acknowledgement{
% The authors appreciate the funding support from Waterloo Institute for Sustainable Energy, Ontario, Canada. Bissan Ghaddar was supported by NSERC Discovery Grant 2017-04185.}
% Acknowledgments here
\ACKNOWLEDGMENT{%
% Enter the text of acknowledgments here
The authors appreciate the funding support from Waterloo Institute for Sustainable Energy, Ontario, Canada. Bissan Ghaddar was supported by NSERC Discovery Grant 2017-04185.
}% Leave this (end of acknowledgment)

\bibliographystyle{apalike}
\addcontentsline{toc}{chapter}{\bibname}
\bibliography{ref}
\end{document}